\documentclass[preprint,11pt]{elsarticle}
\usepackage{makeidx}
\usepackage{amssymb}
\usepackage{amsfonts}
\usepackage{amsmath}

\newtheorem{theorem}{Theorem}
\newtheorem{corollary}[theorem]{Corollary}
\newtheorem{definition}[theorem]{Definition}
\newtheorem{example}[theorem]{Example}
\newtheorem{lemma}[theorem]{Lemma}
\newtheorem{proposition}[theorem]{Proposition}
\newtheorem{remark}[theorem]{Remark}
\newenvironment{proof}[1][Proof]{\noindent\textbf{#1.} }{\ \rule{0.5em}{0.5em}}
\input{tcilatex}
\begin{document}

\begin{frontmatter}



\title{\large On the Rate of Convergence of Weak Euler Approximation for Nondegenerate SDEs Driven by L\'{e}vy Processes}


\author{Remigijus Mikulevi\v{c}ius and Changyong Zhang}

\address{University of Southern California, Los Angeles, USA}

\begin{abstract}
The paper studies the rate of convergence of the weak Euler approximation for solutions to SDEs
driven by L\'{e}vy processes, with H\"{o}lder-continuous coefficients. It investigates the dependence of the rate on the regularity of coefficients and
driving processes. The equation considered has a nondegenerate main part driven by a spherically-symmetric stable process.

\end{abstract}


\end{frontmatter}




\section{\textrm{Introduction}}


The paper studies the weak Euler approximation for solutions to SDEs driven
by L\'{e}vy processes with a nondegenerate main part. The goal is to
investigate the dependence of the convergence rate on the regularity of
coefficients and driving processes. We use the method developed in \cite%
{mikz1} and \cite{mikz2}. For the sake of completeness we repeat some
arguments. A methodical novelty is that contrary to \cite{mikz1} and \cite%
{mikz2} we do not use Fourier transform. Also, the whole H\"{o}lder-Zygmund
scale is covered.

\subsection{Nondegenerate SDEs Driven by L\'{e}vy Processes}

Let $(\Omega ,\mathcal{F},\mathbf{P})$ be a complete probability space with
a filtration $\mathbb{F}=\{\mathcal{F}_{t}\}_{t\in \lbrack 0,T]}$ of $\sigma 
$-algebras satisfying the usual conditions and $\alpha \in (0,2]$ be fixed.
Consider the following model in $\mathbf{R}^{d}$: 
\begin{equation}
X_{t}=X_{0}+\int_{0}^{t}a(X_{s})ds+\int_{0}^{t}b(X_{s-})dU_{s}^{\alpha
}+\int_{0}^{t}G(X_{s-})dZ_{s},t\in \lbrack 0,T],  \label{one}
\end{equation}%
where $a(x)=(a^{i}(x))_{1\leq i\leq d}$, $b(x)=(b^{ij}(x))_{1\leq i,j\leq d}$%
, $G(x)=(G^{ij}(x))_{1\leq i\leq d,1\leq j\leq m}$, $x\in \mathbf{R}^{d}$
are measurable and bounded, with $a=0$ if $\alpha \in (0,1)$ and $b$ being
nondegenerate. The main part of the equation is driven by $U^{\alpha
}=\{U_{t}^{\alpha }\}_{t\in \lbrack 0,T]}$, a standard $d$-dimensional
spherically-symmetric $\alpha $-stable process:%
\begin{equation*}
U_{t}^{\alpha }=\int_{0}^{t}\int (1-\bar{\chi}_{\alpha
}(y))yp_{0}(ds,dy)+\int_{0}^{t}\int \bar{\chi}_{\alpha
}(y)yq_{0}(ds,dy),\alpha \in (0,2),
\end{equation*}%
where $\bar{\chi}_{\alpha }(y)=\mathbf{1}_{\{\alpha \in (1,2)\}}+\mathbf{1}%
_{\{\alpha =1\}}\chi _{\left\{ |y|\leq 1\right\} }$ and $~p_{0}(dt,dy)$ is a
Poisson point measure on $[0,\infty )\times \mathbf{R}_{0}^{d}$ ($\mathbf{R}%
_{0}^{d}=\mathbf{R}^{d}\backslash \left\{ 0\right\} $) with 
\begin{equation*}
\mathbf{E}[p_{0}(dt,dy)]=\frac{dtdy}{|y|^{d+\alpha }},\quad
q_{0}(dt,dy)=p_{0}(dt,dy)-\frac{dtdy}{|y|^{d+\alpha }}.
\end{equation*}%
If $\alpha =2$, $U^{\alpha }$ is the standard Wiener process in $\mathbf{R}%
^{d}$. The last term is driven by $Z=\{Z_{t}\}_{t\in \lbrack 0,T]}$, an $m$%
-dimensional L\'{e}vy process whose characteristic function is $\exp \left\{
t\eta (\xi )\right\} $ with 
\begin{equation*}
\eta (\xi ) = \int_{\mathbf{R}_{0}^{m}}\big[e^{i(\xi ,y)}-1-i(\xi ,y)\chi
_{\left\{ |y|\leq 1\right\} }\mathbf{1}_{\{\alpha \in (1,2]\}}\big]\pi (dy).
\end{equation*}
Hence, 
\begin{equation}
Z_{t}=\int_{0}^{t}\int (1-\chi _{\alpha }(y))yp(ds,dy)+\int_{0}^{t}\int \chi
_{\alpha }(y)yq(ds,dy),  \notag
\end{equation}%
where $\chi _{\alpha }(y)=\mathbf{1}_{\{\alpha \in (1,2]\}}\chi _{\left\{
|y|\leq 1\right\} }$, $p(dt,dy)$ is a Poisson point measure on $[0,\infty
)\times \mathbf{R}_{0}^{m}$ with $\mathbf{E}[p(dt,dy)]=\pi (dy)dt$, and $%
q(dt,dy)=p(dt,dy)-\pi (dy)dt$ is the centered Poisson measure. It is assumed
that 
\begin{equation*}
\int (|y|^{\alpha }\wedge 1)\pi (dy)<\infty .
\end{equation*}

\subsection{Motivation}

The process defined in (\ref{one}) is used as a mathematical model for
random dynamic phenomena in applications arising from fields such as finance
and insurance, to capture continuous and discontinuous uncertainty. For many
applications, the practical computation of functionals of the type $F=%
\mathbf{E}[g(X_{T})]$ and $F=\mathbf{E}\big[\int_{0}^{T}f(X_{s})ds\big]$
plays an important role. For instance in finance, derivative prices can be
expressed by such functionals. However in reality, a stochastic differential
equation does not always have a closed-form solution. In such cases, in
order to evaluate $F$, an alternative option is to numerically approximate
the It\^{o} process $X$ by a discrete-time Monte-Carlo simulation, which has
been widely applied. The simplest and the most commonly-used scheme is the
weak Euler approximation.

Let the time discretization $\{\tau _{i},i=0,\ldots ,n_{T}\}$ of the
interval $[0,T]$ with maximum step size $\delta \in (0,1)$ be a partition of 
$[0,T]$ such that $0=\tau _{0}<\tau _{1}<\dots <\tau _{n_{T}}=T$ and $%
\max_{i}(\tau _{i}-\tau _{i-1})\leq \delta .$ \ The Euler approximation of $%
X $ is an $\mathbb{F}$-adapted stochastic process $Y=\{Y_{t}\}_{t\in \lbrack
0,T]}$ defined by the stochastic equation%
\begin{equation}
Y_{t}=X_{0}+\int_{0}^{t}a(Y_{\tau _{i_{s}}})ds+\int_{0}^{t}b(Y_{\tau
_{i_{s}}})dU_{s}^{\alpha }+\int_{0}^{t}G(Y_{\tau _{i_{s}}})dZ_{s},t\in
\lbrack 0,T],  \label{du}
\end{equation}%
where $\tau _{i_{s}}=\tau _{i}$\ if $s\in \lbrack \tau _{i},\tau
_{i+1}),i=0,\ldots ,n_{T}-1.$ Contrary to those in (\ref{one}), the
coefficients in (\ref{du}) are piecewise constants in each time interval of $%
[\tau _{i},\tau _{i+1}).$

The weak Euler approximation $Y$\ is said to converge with order $\kappa >0$%
\ if for each bounded smooth function $g$ with bounded derivatives, there
exists a constant $C$, depending only on $g$, such that 
\begin{equation*}
|\mathbf{E}[g(Y_{T})] - \mathbf{E}[g(X_{T})]|\leq C\delta ^{\kappa },
\end{equation*}%
where $\delta >0$\ is the maximum step size of the time discretization.

In the literature, the weak Euler approximation of stochastic differential
equations with smooth coefficients has been consistently studied. For
diffusion processes ($\alpha =2)$, Milstein was one of the first to
investigate the order of weak convergence and derived $\kappa =1$~\cite%
{Mil79, Mil86}. Talay considered a class of the second order approximations
for diffusion processes~\cite{Tal84, Tal86}. For It\^{o} processes with jump
components, Mikulevi\v{c}ius \& Platen showed the first-order convergence in
the case in which the coefficient functions possess fourth-order continuous
derivatives~\cite{MiP88}. Platen and Kloeden \& Platen studied not only
Euler but also higher order approximations~\cite{KlP00, Pla99}. Protter \&
Talay analyzed the weak Euler approximation for 
\begin{equation}
X_{t}=X_{0}+\int_{0}^{t}G(X_{s-})dZ_{s},t\in \lbrack 0,T],  \label{pt}
\end{equation}%
where $Z_{t}=(Z_{t}^{1},\ldots ,Z_{t}^{m})$ is a L\'{e}vy process and $%
G=(G^{ij})_{1\leq i\leq d,1\leq j\leq m}$ is a measurable and bounded
function~\cite{PrT97}. They showed the order of convergence $\kappa =1,$
provided that $G$ and $g$\ are smooth and the L\'{e}vy measure of $Z$\ has
finite moments of sufficiently high order. Because of this, the main
theorems in \cite{PrT97} do not apply to (\ref{one}). On the other hand, (%
\ref{one}) with a nondegenerate matrix $b$ does not cover (\ref{pt}), which
can degenerate completely.

In general, the coefficients and the test function $g$\ do not always have
the smoothness properties assumed in the papers cited above. Mikulevi\v{c}%
ius \& Platen proved that there still exists some order of convergence of
the weak Euler approximation for nondegenerate diffusion processes under H%
\"{o}lder conditions on the coefficients and $g$~\cite{MiP911}. Kubilius \&
Platen and Platen \& Bruti-Liberati considered a weak Euler approximation in
the case of a nondegenerate diffusion process with a finite number of jumps
in finite time intervals~\cite{KuP01, PlB10}.

In this paper, we investigate the dependence of the rate of convergence on
the H\"{o}lder regularity of coefficients and the driving processes. For a
driving process, the variation of the process can be regarded as a part of
its regularity. In this sense, Wiener process is the worse, most
\textquotedblleft chaotic\textquotedblright , among $\alpha $-stable
processes. Also, as pointed out in \cite{PrT97}, the tails of L\'{e}vy
processes influence the convergence rate as well.

\subsection{Examples}

For $\beta >0\,,$ denote $C^{\beta }(\mathbf{R}^{d})$ the H\"{o}lder-Zygmund
space, and $\tilde{C}^{\beta }(\mathbf{R}^{d})$ the Lipshitz space ($\tilde{C%
}^{\beta }(\mathbf{R}^{d})=C^{\beta }(\mathbf{R}^{d})$ if $\beta \notin 
\mathbf{N}$, see Section~\ref{sec:operator} for definitions). Let us look at
two examples.

\begin{example}
\label{ex1}$($see Corollary \ref{co1}$)$ Assume $\beta \leq \alpha $, the
coefficients $a^{i},b^{ij}\in \tilde{C}^{\beta }(\mathbf{R}^{d})$, $%
G^{ij}\in \tilde{C}^{\frac{\beta }{\alpha \wedge 1}}(\mathbf{R}^{d})$, $%
\inf_{x}|\det b(x)|>0,$ and 
\begin{equation*}
\int_{\mathbf{R}^{m}}|y|^{\alpha }\pi (dy)<\infty ,
\end{equation*}%
where $\pi $ is the L\'{e}vy measure of the driving process $Z$. Then it
holds that%
\begin{eqnarray*}
|\mathbf{E}[g(Y_{T})]-\mathbf{E}[g(X_{T})]| &\leq &C|g|_{\alpha +\beta
}r(\delta ,\alpha ,\beta ), \\
|\mathbf{E}\big[\int_{0}^{T}f(Y_{\tau _{i_{s}}})ds\big]-\mathbf{E}\big[%
\int_{0}^{T}f(X_{s})ds\big]| &\leq &C|f|_{\beta }r(\delta ,\alpha ,\beta ),
\end{eqnarray*}%
where 
\begin{equation*}
r(\delta ,\alpha ,\beta )=\left\{ 
\begin{array}{cc}
\delta ^{\frac{\beta }{\alpha }} & \text{if }\beta <\alpha , \\ 
\delta (1+\left\vert \ln \delta \right\vert ) & \text{if }\beta =\alpha%
\end{array}%
\right.
\end{equation*}
\end{example}

\begin{example}
\label{ex2}$($see Corollary \ref{co2}$)$ Consider the jump-diffusion case $%
(\alpha =2)$ 
\begin{equation*}
X_{t}=X_{0}+\int_{0}^{t}a(X_{s})ds+\int_{0}^{t}b(X_{s})dW_{s}+%
\int_{0}^{t}G(X_{s-})dZ_{s},t\in \lbrack 0,T],
\end{equation*}%
where $W=\{W_{t}\}_{t\in \lbrack 0,T]}$ is a standard Wiener process. Assume 
$a,b^{ij}\in \tilde{C}^{\beta }(\mathbf{R}^{d})$, $\inf_{x}|\det b(x)|>0$,
and there exists $\mu \in (0,3)$ such that%
\begin{equation*}
\int_{|y|\leq 1}|y|^{2}\pi (dy)+\int_{|y|>1}|y|^{\mu }\pi (dy)<\infty 
\end{equation*}%
Let $G^{ij}\in \tilde{C}^{\frac{\beta }{\mu \wedge 1}}(\mathbf{R}^{d})$.
Then it holds that%
\begin{eqnarray*}
|\mathbf{E}[g(Y_{T})]-\mathbf{E}[g(X_{T})]| &\leq &C|g|_{\alpha +\beta
}r(\delta ), \\
|\mathbf{E}\big[\int_{0}^{T}f(Y_{\tau _{i_{s}}})ds\big]-\mathbf{E}\big[%
\int_{0}^{T}f(X_{s})ds\big]| &\leq &C|f|_{\beta }r(\delta ),
\end{eqnarray*}%
where%
\begin{equation*}
r(\delta )=\left\{ 
\begin{array}{cc}
\delta ^{\frac{\beta \wedge \mu }{2}} & \text{if }\mu <2, \\ 
\delta (1+|\ln \delta |) & \text{if }\mu =\beta =2, \\ 
\delta  & \text{if }\mu >2,\beta >2%
\end{array}%
\right. 
\end{equation*}%
The assumption $G^{ij}\in \tilde{C}^{\frac{\beta }{\mu \wedge 1}}(\mathbf{R}%
^{d})$ shows that if $\mu <1$, the heavy tail of $\pi $ can be balanced by a
higher regularity of $G^{ij}$.
\end{example}

As in \cite{MiP911}, this paper employs the idea of Talay (see \cite{Tal84})
and uses the solution to the backward Kolmogorov equation associated with $%
X_{t},$ It\^{o}'s formula, and one-step estimates (see Section~\ref%
{sec:outline} for the outline of the proof).

The paper is organized as follows. In Section 2, the main result is stated
and the proof is outlined. In Section 3, we present the essential technical
results, followed by the proof of the main theorem in Section 4.


\section{Notation and Main Result}



\subsection{Main Result and Notation}

The main result of this paper is the following statement.

\begin{theorem}
\label{thm:main}Let $\beta \in (0,3)$, $0<\beta \leq \mu <\alpha +\beta $
and 
\begin{equation*}
\int_{|y|\leq 1}|y|^{\alpha }\pi (dy)+\int_{|y|>1}|y|^{\mu }\pi (dy)<\infty .
\end{equation*}%
Assume $\inf_{x}|\det b(x)|>0$ and $a^{i},b^{ij}\in \tilde{C}^{\beta }(%
\mathbf{R}^{d}),G^{ij}\in \tilde{C}^{\frac{\beta }{\mu \wedge 1}}(\mathbf{R}%
^{d})$. Then there exists a constant $C$ such that for all $g\in C^{\alpha
+\beta }(\mathbf{R}^{d})$, $f\in C^{\beta }(\mathbf{R}^{d})$,%
\begin{eqnarray*}
|\mathbf{E}[g(Y_{T})]-\mathbf{E}[g(X_{T})]| &\leq &C|g|_{\alpha +\beta
}r(\delta ,\alpha ,\beta ), \\
|\mathbf{E}\big[\int_{0}^{T}f(Y_{\tau _{i_{s}}})ds\big]-\mathbf{E}\big[%
\int_{0}^{T}f(X_{s})ds\big]| &\leq &C|f|_{\beta }r(\delta ,\alpha ,\beta ),
\end{eqnarray*}%
where%
\begin{equation*}
r(\delta ,\alpha ,\beta )=\left\{ 
\begin{array}{cl}
\delta ^{\frac{\beta }{\alpha }}, & \beta <\alpha , \\ 
\delta (1+\left\vert \ln \delta \right\vert ), & \beta =\alpha , \\ 
\delta , & \beta >\alpha .%
\end{array}%
\right.
\end{equation*}
\end{theorem}

Applying Theorem~\ref{thm:main} to the case $\alpha =\mu $ and the case of
heavier tails results in Corollary~\ref{co1} and Corollary~\ref{co2},
respectively.

\begin{corollary}
\label{co1}Let $\beta \in (0,3)$, $\beta \leq \alpha $, and 
\begin{equation*}
\int |y|^{\alpha }\pi (dy)<\infty .
\end{equation*}%
Assume $a^{i},b^{ij}\in \tilde{C}^{\beta }(\mathbf{R}^{d})$, $G^{ij}\in 
\tilde{C}^{\frac{\beta }{\alpha \wedge 1}}(\mathbf{R}^{d})$, and $%
\inf_{x}|\det b(x)|>0.$ Then there exists a constant $C$ such that for all $%
g\in C^{\alpha +\beta }(\mathbf{R}^{d})$, $f\in C^{\beta }(\mathbf{R}^{d})$,%
\begin{eqnarray*}
|\mathbf{E}[g(Y_{T})]-\mathbf{E}[g(X_{T})]| &\leq &C|g|_{\alpha +\beta
}r(\delta ,\alpha ,\beta ), \\
|\mathbf{E}\big[\int_{0}^{T}f(Y_{\tau _{i_{s}}})ds\big]-\mathbf{E}\big[%
\int_{0}^{T}f(X_{s})ds\big]| &\leq &C|f|_{\beta }r(\delta ,\alpha ,\beta ),
\end{eqnarray*}%
where 
\begin{equation*}
r(\delta ,\alpha ,\beta )=\left\{ 
\begin{array}{cc}
\delta ^{\frac{\beta }{\alpha }} & \text{if }\beta <\alpha , \\ 
\delta (1+\left\vert \ln \delta \right\vert ) & \text{if }\beta =\alpha%
\end{array}%
\right.
\end{equation*}
\end{corollary}

\begin{corollary}
\label{co2}Let $\beta \in (0,3)$, $0<\beta \leq \mu <\alpha $, and 
\begin{equation*}
\int_{|y|\leq 1}|y|^{\alpha }\pi (dy)+\int_{|y|>1}|y|^{\mu }\pi (dy)<\infty .
\end{equation*}%
Let $\inf_{x}|\det b(x)|>0,a^{i},b^{ij}\in \tilde{C}^{\beta }(\mathbf{R}%
^{d}) $ and $G^{ij}\in \tilde{C}^{\frac{\beta }{\mu \wedge 1}}(\mathbf{R}%
^{d})$. Then there exists a constant $C$ such that for all $g\in C^{\alpha
+\beta }(\mathbf{R}^{d})$, $f\in C^{\beta }(\mathbf{R}^{d})$,%
\begin{eqnarray*}
|\mathbf{E}[g(Y_{T})]-\mathbf{E}[g(X_{T})]| &\leq &C|g|_{\alpha +\beta
}\delta ^{\frac{\beta \wedge \mu }{\alpha }}, \\
|\mathbf{E}\big[\int_{0}^{T}f(Y_{\tau _{i_{s}}})ds\big]-\mathbf{E}\big[%
\int_{0}^{T}f(X_{s})ds\big]| &\leq &C|f|_{\beta }\delta ^{^{\frac{\beta
\wedge \mu }{\alpha }}}.
\end{eqnarray*}
\end{corollary}

Denote $H=[0,T]\times \mathbf{R}^{d}$, $\mathbf{N}=\{0,1,2,\ldots \}$, $%
\mathbf{R}_{0}^{d}=\mathbf{R}^{d}\backslash \{0\}$. For $x,y\in \mathbf{R}%
^{d}$, write $(x,y)=\sum_{i=1}^{d}x_{i}y_{i}$. For $(t,x)\in H,$ multiindex $%
\gamma \in \mathbf{N}^{d}$ with $D^{\gamma }=\frac{\partial ^{|\gamma |}}{%
\partial x_{1}^{\gamma _{1}}\ldots \partial x_{d}^{\gamma _{d}}}$, and $%
i,j=1,\ldots ,d$, denote 
\begin{eqnarray*}
\partial _{t}u(t,x) &=&\frac{\partial }{\partial t}u(t,x), \ D^{k}u(t,x)=%
\big(D^{\gamma }u(t,x)\big)_{|\gamma |=k},k\in \mathbf{N}\text{,} \\
\partial _{i}u(t,x) &=&u_{x_{i}}(t,x)=\frac{\partial }{\partial x_{i}}%
u(t,x), \ \partial _{ij}^{2}u(t,x)=u_{x_{i}x_{j}}(t,x)=\frac{\partial ^{2}}{%
\partial x_{i}x_{j}}u(t,x), \\
\partial _{x}u(t,x) &=&\nabla u(t,x)=\big(\partial _{1}u(t,x),\dots
,\partial _{d}u(t,x)\big), \\
\partial ^{2}u(t,x) &=&\Delta u(t,x)=\sum_{i=1}^{d}\partial _{ii}^{2}u(t,x).
\end{eqnarray*}

$C=C(\cdot ,\ldots ,\cdot )$ denotes constants depending only on quantities
appearing in parentheses. In a given context the same letter is (generally)
used to denote different constants depending on the same set of arguments.

\subsection{Outline of Proof}

\label{sec:outline}

Due to the lack of regularity, standard techniques such as stochastic flows
cannot be applied to prove Theorem~\ref{thm:main}. Instead, as in \cite%
{MiP911}, the solution to the backward Kolmogorov equation associated with $%
X_{t}$ is used. In the following, the operators of the Kolmogorov equation
associated with $X_{t}$ are first defined.

For $u\in C^{\alpha +\beta }(H)$, denote 
\begin{eqnarray*}
A_{z}u(t,x) &=&\mathbf{1}_{\{\alpha =1\}}(a(z),\nabla _{x}u(t,x))+\mathbf{1}%
_{\{\alpha =2\}}\frac{1}{2}\sum_{i,j=1}^{d}D^{ij}(z)\partial _{ij}^{2}u(t,x)
\\
&&+\mathbf{1}_{\{\alpha \in (0,2)\}}\int [u(t,x+b(z)y)-u(t,x)-(\nabla
u(t,x),b(z)y)\chi _{\alpha }(y)]\frac{dy}{|y|^{d+\alpha }}, \\
Au(t,x) &=&A_{x}u(t,x)=A_{z}u(t,x)|_{z=x},
\end{eqnarray*}%
with $\chi _{\alpha }(y)=\mathbf{1}_{\{\alpha \in (1,2)\}}+\mathbf{1}%
_{\{\alpha =1\}}\chi _{\left\{ |y|\leq 1\right\} },D=b^{\ast }b,$ and%
\begin{eqnarray*}
B_{z}u(t,x) &=&\mathbf{1}_{\{\alpha \in (1,2]\}}(a(z),\nabla
_{x}u(t,x))+\int_{\mathbf{R}_{0}^{m}}\big[u(t,x+G(z)y))-u(t,x) \\
&&-\mathbf{1}_{\{\alpha \in (1,2]\}}\mathbf{1}_{\left\{ |y|\leq 1\right\}
}(\nabla _{x}u(t,x),G(z)y))\big]\pi (dy), \\
Bu(t,x) &=&B_{x}u(t,x)=B_{z}u(t,x)|_{z=x}.
\end{eqnarray*}

Applying It\^{o}'s formula to $X_{t}$ and $u\in C_{0}^{\infty }(\mathbf{R}%
^{d})$, we find that 
\begin{equation*}
u(X_{t})-\int_{0}^{t}Au(X_{s})ds-\int_{0}^{t}Bu(X_{s})ds,t\in \lbrack 0,T]
\end{equation*}%
is a martingale.

\begin{remark}
\label{lrenew1}More precisely, under assumptions of Theorem \ref{thm:main},
there exists a unique weak solution to equation \textup{(\ref{one})} and the
stochastic process 
\begin{equation*}
u(X_{t})-\int_{0}^{t}(A+B)u(X_{s})ds,\forall u\in C^{\alpha +\beta }(\mathbf{%
R}^{d})  \label{ff25}
\end{equation*}%
is a martingale~\textup{\cite{MiP923}}. The operator $\mathcal{L}=A+B$ is
the generator of $X_{t}$ defined in \textup{(\ref{one})}; $A$ is the
principal part of $\mathcal{L}$ and $B$ is the lower order or subordinated
part of $\mathcal{L}$.
\end{remark}

If $v(t,x),(t,x)\in H$ satisfies the backward Kolmogorov equation 
\begin{eqnarray*}
\big(\partial _{t}+A+B\big)v(t,x) &=&0,\quad 0\leq t\leq T, \\
v(T,x) &=&g(x),
\end{eqnarray*}%
then as interpreted in Section 4, by It\^{o}'s formula 
\begin{equation*}
\mathbf{E}[g(Y_{T})]-\mathbf{E}[g(X_{T})]=\mathbf{E}[v(T,Y_{T})-v(0,Y_{0})]=%
\mathbf{E}\big[\int_{0}^{T}(\partial _{t}+\mathcal{L}_{Y_{\tau
_{i_{s}}}})v(s,Y_{s})ds\big].
\end{equation*}%
The regularity of $v$ determines the one-step estimate and the rate of
convergence of the approximation. For $\beta \in (0,1)$, the results for the
Kolmogorov equation in H\"{o}lder classes are available~\cite{MiP922, MiP09}%
. In a standard way the results can be extended to the case $\beta >1$. The
main difficulty is to derive the one-step estimates (see Lemma~\ref%
{lem:expect}).


\section{Backward Kolmogorov Equation}


In H\"{o}lder-Zygmund spaces, consider the backward Kolmogorov equation
associated with $X_{t}$: 
\begin{eqnarray}
\big(\partial _{t}+A+B\big)v(t,x) &=&f(t,x),  \label{eq1} \\
v(T,x) &=&0.  \notag
\end{eqnarray}%
The regularity of its solution is essential for the one-step estimate which
determines the rate of convergence.

\begin{definition}
\label{def1}Let $f$ be a measurable and bounded function on $\mathbf{R}^{d}$%
. We say that $u\in C^{\alpha +\beta }(H)$ is a solution to $(\ref{eq1})$ if 
\begin{equation}
u(t,x)=\int_{t}^{T}\big[\mathcal{L}u(s,x)-f(s,x)\big]ds,\forall (t,x)\in H.
\label{defs}
\end{equation}
\end{definition}


The following theorem is the main result of this section.

\begin{theorem}
\label{thm:StoCP}Let $\beta \in (0,3)$, $0<\beta \leq \mu <\alpha +\beta $,
and 
\begin{equation*}
\int_{|y|\leq 1}|y|^{\alpha }\pi (dy)+\int_{|y|>1}|y|^{\mu }\pi (dy)<\infty .
\end{equation*}%
Assume $a^{i},b^{ij}\in \tilde{C}^{\beta }(\mathbf{R}^{d}),G^{ij}\in \tilde{C%
}^{\frac{\beta }{\mu \wedge 1}}(\mathbf{R}^{d}),\inf_{x}|\det b(x)|>0$. Then
for each $f\in C^{\beta }(\mathbf{R}^{d})$, there exists a unique solution $%
v\in C^{\alpha +\beta }(H)$ to \textup{(\ref{eq1})} and a constant $C$
independent of $f$ such that $|u|_{\alpha +\beta }\leq C|f|_{\beta }.$
\end{theorem}

An immediate consequence of Theorem~\ref{thm:StoCP} is the following
statement.

\begin{corollary}
\label{lcornew1}Let $\beta \in (0,3)$, $0<\beta \leq \mu <\alpha +\beta $,
and 
\begin{equation*}
\int_{|y|\leq 1}|y|^{\alpha }\pi (dy)+\int_{|y|>1}|y|^{\mu }\pi (dy)<\infty .
\end{equation*}%
Let $a^{i},b^{ij}\in \tilde{C}^{\beta }(\mathbf{R}^{d}),G^{ij}\in \tilde{C}^{%
\frac{\beta }{\mu \wedge 1}}(\mathbf{R}^{d}),\inf_{x}|\det b(x)|>0$. Then
for each $f\in C^{\beta }(\mathbf{R}^{d})$ and $g\in C^{\alpha +\beta }(%
\mathbf{R}^{d})$, there exists a unique solution $v\in C^{\alpha +\beta }(H)$
to the Cauchy problem 
\begin{eqnarray}
\big(\partial _{t}+A+B\big)v(t,x) &=&f(x),  \label{maf8} \\
v(T,x) &=&g(x)  \notag
\end{eqnarray}%
and $|v|_{\alpha +\beta }\leq C(|f|_{\beta }+|g|_{\alpha +\beta })$ with a
constant $C$ independent of $f$ and $g$.
\end{corollary}

To prove Theorem \ref{thm:StoCP} and Corollary \ref{lcornew1}, in a standard
way the equation with constant coefficients is first solved. Then variable
coefficients are handled by using partition of unity and deriving apriori
Schauder estimates in H\"{o}lder-Zygmund spaces. Finally, the continuation
by parameter method is applied to extend solvability of an equation with
constant coefficients to (\ref{eq1}).

\subsection{Kolmogorov Equation with Constant Coefficients}

It is convenient to rewrite the principal operator $A$ by changing the
variable of integration in the integral part:%
\begin{eqnarray*}
A_{z}u(t,x) &=&\mathbf{1}_{\{\alpha =1\}}(a(z),\nabla _{x}u(t,x))+\mathbf{1}%
_{\{\alpha =2\}}\sum_{i,j=1}^{d}D^{ij}(z)\partial _{ij}^{2}u(t,x) \\
&&+\mathbf{1}_{\{\alpha \in (0,2)\}}\int [u(t,x+y)-u(t,x)-(\nabla
u(t,x),y)\chi _{\alpha }(y)]m(z,y)\frac{dy}{|y|^{d+\alpha }},
\end{eqnarray*}%
where $D=b^{\ast }b,$%
\begin{equation}
m(z,y)=\frac{1}{|\det b(z)|}\frac{1}{|b(z)^{-1}\frac{y}{|y|}|^{d+\alpha }}%
,\alpha \in (0,2).  \label{foo}
\end{equation}%
Obviously, 
\begin{equation}
\int_{S^{d-1}}ym(\cdot ,y)\mu _{d-1}(dy)=0.  \label{ff27}
\end{equation}%
Here $S^{d-1}$ is the unit sphere in $\mathbf{R}^{d}$ and $\mu _{d-1}$ is
the Lebesgue measure.

For $z_{0}\in \mathbf{R}^{d}$, denote $A^{0}u(x)=A_{z_{0}}u(x)$. Consider a
backward Kolmogorov equation with constant coefficients and $\lambda \geq 0,$%
\begin{eqnarray}
\big(\partial _{t}+A^{0}-\lambda \big)v(t,x) &=&f(x),  \label{eq0} \\
v(T,x) &=&0.  \notag
\end{eqnarray}

\begin{proposition}
\label{prop1}Let $\beta >0$ and $f\in C^{\beta }(\mathbf{R}^{d})$. Assume
there are constants $c_{1},K>0$ such that for all $z\in \mathbf{R}^{d},$ 
\begin{equation*}
|\det b(z)|\geq c_{1},\quad \mathbf{1}_{\{\alpha =1\}}|a(z)|+|b(z)|\leq K.
\end{equation*}%
Then there exists a unique solution $u\in C^{\alpha +\beta }(H)$ to (\ref%
{eq0}) and 
\begin{equation}
|u|_{\alpha +\beta }\leqslant C|f|_{\beta },  \label{1}
\end{equation}%
where the constant $C$ depends only on $\alpha ,\ \beta ,\ T,\ d$, $c_{1},\
K.$ Moreover,%
\begin{equation}
|u|_{\beta }\leq C(\alpha ,d)(\lambda ^{-1}\wedge T)|f|_{\beta }  \label{2}
\end{equation}%
and there exists a constant $C$ such that for all $s\leq t\leq T,$ 
\begin{equation}
|u(t,\cdot )-u(s,\cdot )|_{\frac{\alpha }{2}+\beta }\leq C(t-s)^{\frac{1}{2}%
}|f|_{\beta }.  \label{3}
\end{equation}
\end{proposition}

To derive Proposition~\ref{prop1}, some auxiliary results are presented
first.

\subsubsection{Continuity of Operator $A^{0}$ in H\"{o}lder-Zygmund Spaces}

\label{sec:operator}

To show that operator $A^{0}$ is continuous in H\"{o}lder-Zygmund spaces $%
C^{\beta }(\mathbf{R}^{d})$, first recall their definition.

For $\beta =[\beta ]^{-}+\left\{ \beta \right\} ^{+}>0$, where $[\beta
]^{-}\in \mathbf{N}$ and $\left\{ \beta \right\} ^{+}\in (0,1]$, let $%
C^{\beta }(H)$ denote the space of measurable functions $u$ on $H$ such that
the norm 
\begin{eqnarray*}
|u|_{\beta } &=&\sum_{|\gamma |\leq \lbrack \beta ]^{-}}\sup_{(t,x)\in
H}|D_{x}^{\gamma }u(t,x)|+\mathbf{1}_{\{\{\beta \}^{+}<1\}}\sup_{\substack{ %
|\gamma |=[\beta ]^{-},  \\ t,x\neq \tilde{x}}}\frac{|D_{x}^{\gamma
}u(t,x)-D_{x}^{\gamma }u(t,\tilde{x})|}{|x-\tilde{x}|^{\{\beta \}^{+}}} \\
&&+\mathbf{1}_{\{\{\beta \}^{+}=1\}}\sup_{\substack{ |\gamma |=[\beta ]^{-}, 
\\ t,x,h\neq 0}}\frac{|D_{x}^{\gamma }u(t,x+h)+D_{x}^{\gamma
}u(t,x-h)-2D_{x}^{\gamma }u(t,x)|}{|h|^{\{\beta \}^{+}}}
\end{eqnarray*}%
is finite. Accordingly, $C^{\beta }(\mathbf{R}^{d})$ denotes the
corresponding space of functions on $\mathbf{R}^{d}$.~The classes $C^{\beta
} $ coincide with H\"{o}lder spaces if $\beta \notin \mathbf{N}$ (see 1.2.2
of \cite{Tri92}).

For $v\in C^{\beta }(\mathbf{R}^{d})$ with $\beta \in (0,1]$, denote%
\begin{eqnarray*}
|v|_{0} &=&\sup_{x}|v(x)|,\  \\
\left[ v\right] _{\beta } &=&\sup_{x\neq y}\frac{|v(x)-v(y)|}{|x-y|^{\beta }}%
\text{ if }\beta \in (0,1), \\
\left[ v\right] _{\beta } &=&\sup_{x,h\neq 0}\frac{|v(x+h)-2v(x)+v(x-h)|}{%
|x-y|^{\beta }}\text{ if }\beta =1.
\end{eqnarray*}%
Similarly we define the spaces $\tilde{C}^{\beta }(\mathbf{R}^{d})$: $\tilde{%
C}^{\beta }(\mathbf{R}^{d})=C^{\beta }(\mathbf{R}^{d})$ if $\beta >0,\beta
\notin \mathbf{N}$, and $\tilde{C}^{k}(\mathbf{R}^{d})$ is the space of all
functions $f$ on $\mathbf{R}^{d}$ having $k-1$ continuous bounded
derivatives and such that $D^{\gamma }f,|\gamma |=k-1,$ are Lipshitz. We
introduce the norms in $\tilde{C}^{\beta }(\mathbf{R}^{d})$%
\begin{eqnarray*}
||f||_{\beta } &=&|f|_{\beta }\text{ if }\beta >0,\beta \notin \mathbf{N}, \\
||f||_{k} &=&\sum_{|\gamma |\leq k-1}|D^{\gamma }f|_{0}+\sup_{x\neq
y,|\gamma |=k-1}\frac{|D^{\gamma }f(x)-D^{\gamma }f(y)|}{|x-y|}.
\end{eqnarray*}

For $\alpha \in (0,2)$, define for $v\in C^{\alpha +\beta }(\mathbf{R}^{d})$
the fractional Laplacian%
\begin{equation}
\partial ^{\alpha }v(x)=\int [v(x+y)-v(x)-\left( \nabla v(x),y\right) \chi
_{\alpha }(y)]\frac{dy}{|y|^{d+\alpha }},x\in \mathbf{R}^{d}.  \label{fo4}
\end{equation}

For various estimates, the following representation of the difference is
useful.

\begin{lemma}
\label{r1}$($Lemma 2.1 in \textup{\cite{Kom84}}$)$ For $\delta \in (0,1)$
and $u\in C_{0}^{\infty }(\mathbf{R}^{d})$, 
\begin{equation*}
u\big(x+y\big)-u(x)=K\int k^{(\delta )}(y,z)\partial ^{\delta }u(x-z)dz,
\label{eqn:diff_bound}
\end{equation*}%
where $K=K(\delta ,d)$ is a constant, 
\begin{equation*}
k^{(\delta )}(y,z)=|z+y|^{-d+\delta }-|z|^{-d+\delta },
\end{equation*}%
and there exists a constant $C$ such that 
\begin{equation*}
\int |k^{(\delta )}(y,z)|dz\leq C|y|^{\delta },\forall y\in \mathbf{R}^{d}.
\end{equation*}
\end{lemma}

By taking pointwise limit ($\partial ^{\delta}$ is defined by (\ref{fo4}))
and applying the dominated convergence theorem, the statement can be
extended to $u\in C^{\delta }(\mathbf{R}^{d})$.

Let $m(y)$ be a measurable and bounded function on $\mathbf{R}^{d}.$ Define 
\begin{equation*}
L^{m}u(x)=\int_{\mathbf{R}^{d}}\left[ u(x+y)-u(x)-(\nabla u(x), y)\chi
_{\alpha}(y)\right] m(y)\frac{dy}{|y|^{d+\alpha }},u\in C^{\alpha +\beta }.
\end{equation*}

The following statement is proved in~\cite{MiP11} for $\beta \in (0,1]$.

\begin{lemma}
\label{prop2}Let $\alpha \in (0,2)$, $\beta >0$, $u\in C^{\alpha +\beta }(%
\mathbf{R}^{d})$, and $|m|\leq K$. Assume if $\alpha =1$, 
\begin{equation}
\int_{r<|y|\leq 1}ym(y)\frac{dy}{|y|^{d+1}}=0,\forall r\in (0,1).  \label{8}
\end{equation}%
Then there exists a constant $C$ independent of $u$ such that%
\begin{equation*}
|L^{m}u|_{\beta }\leq CK|u|_{\alpha +\beta }.
\end{equation*}
\end{lemma}

\begin{proof}
The result holds for $\beta \in (0,1]$ according to Proposition 11 in \cite%
{MiP11}. If $\beta >1$, and $u\in C^{\alpha +\beta }(\mathbf{R}^{d})$, then
for any multiindex $|\gamma |=[\beta ]$, $D^{\gamma }u\in C^{\alpha +\beta
-[\beta ]}$, and%
\begin{equation*}
|D^{\gamma }\left( L^{m}u\right) |_{\beta -[\beta ]^{-}}=|L^{m}\left(
D^{\gamma }u\right) |_{\beta -[\beta ]^{-}}\leq CK|D^{\gamma }u|_{\alpha
+\beta -[\beta ]^{-}}.
\end{equation*}%
The statement follows.
\end{proof}

\subsubsection{Proof of Proposition \protect\ref{prop1}}

The statement is proved by induction. Given $\alpha \in (0,2]$ and $f\in
C^{\beta }(H)$, for $\beta \in (0,1]$, there exists a unique solution $u\in
C^{\alpha +\beta }(H)$ to the Kolmogorov equation (\ref{eq0}) such that (\ref%
{1})-(\ref{3}) hold~\cite{MiP11}.

Assume the result holds for $\beta \in \bigcup\limits_{l=0}^{n-1}(l,l+1]$, $%
n\in \mathbf{N}$. Let $\beta \in (n,n+1]$, $\tilde{\beta}=\beta -1$, and $%
f\in C^{\beta }(H)$. Then $\tilde{\beta}\in (n-1,n]$, $f\in C^{\tilde{\beta}%
}(H)$, and there exists a unique solution $v\in C^{\alpha +\tilde{\beta}}(H)$%
, $\alpha \in (0,2]$ to the Cauchy problem such that (\ref{1})-(\ref{3})
hold for $v$ with $\tilde{\beta}$. For $h\in \mathbf{R}$ and $k=1,\dots ,d$,
denote 
\begin{equation*}
v_{k}^{h}(t,x)=\frac{v(t,x+he_{k})-v(t,x)}{h},
\end{equation*}%
where $\{e_{k},k=1,\ldots ,d\}$ is the canonical basis in $\mathbf{R}^{d}$.
Obviously, $v_{k}^{h}\in C^{\alpha +\tilde{\beta}}(H)$ and 
\begin{eqnarray}
\big(\partial _{t}+A^{0}-\lambda \big)v_{k}^{h}(t,x) &=&f_{k}^{h}(x),x\in 
\mathbf{R}^{d},  \label{eqn:cauchy} \\
v_{k}^{h}(T,x) &=&0.  \notag
\end{eqnarray}

Since $f\in C^{\beta }(H)$ and 
\begin{equation*}
f_{k}^{h}(t,x)=\int_{0}^{1}\partial _{k}f(t,x+he_{k}s)ds,\forall h\neq 0,
\end{equation*}%
then 
\begin{equation}
|f_{k}^{h}|_{\tilde{\beta}}\leq C|\nabla f|_{\beta -1}\leq C|f|_{\beta }
\label{ff22}
\end{equation}%
with a constant $C$ independent of $h$. Since $v\in C^{\alpha +\tilde{\beta}%
}(H)$, then $v_{k}^{h}\in C^{\alpha +\tilde{\beta}}(H)$. By (\ref{ff22}) and
the induction assumption, the estimates (\ref{1})-(\ref{3}) hold for $%
v_{k}^{h}$ with a constant independent of $h$. Hence $v_{k}^{h}(t,x)$ are
equicontinuous in $(t,x)$. By the Arzel\`{a}-Ascoli theorem, for each $%
h_{n}\rightarrow 0$, there exist a subsequence $\{h_{n_{j}}\}$ and
continuous functions $v_{k}(t,x),(t,x)\in H,k=1,\ldots ,d$, such that $%
v_{k}^{h_{n_{j}}}(t,x)\rightarrow v_{k}(t,x)$ uniformly on compact subsets
of $H$ as $j\rightarrow \infty $. Therefore, $v_{k}\in C^{\alpha +\tilde{%
\beta}}$ and $|v_{k}|_{\alpha +\tilde{\beta}}\leq C|f|_{\beta },k=1,\ldots
,d.$

It then follows from passing to the limit in the integral form of (\ref%
{eqn:cauchy}) (see (\ref{defs})) and the dominated convergence theorem that $%
u_{k}$ is the unique solution to 
\begin{eqnarray*}
\big(\partial _{t}+A^{0}-\lambda \big)v_{k}(t,x) &=&\partial _{k}f(t,x), \\
v_{k}(T,x) &=&0,k=1,\dots ,d
\end{eqnarray*}%
and so $v_{k}^{h_{n}}(t,x)\rightarrow v_{k}(t,x),\forall h_{n}\rightarrow 0.$
Hence, 
\begin{equation*}
v_{k}(t,x)=\lim_{h\rightarrow 0}v_{k}^{h}(t,x)=\lim_{h\rightarrow 0}\frac{%
v(t,x+he_{k})-v(t,x)}{h}=\partial _{k}v(t,x),
\end{equation*}%
$\partial _{k}v\in C^{\alpha +\tilde{\beta}}(H),k=1,\dots ,d$, and $|\nabla
v|_{\alpha +\tilde{\beta}}\leq C|f|_{\beta }.$ Therefore, $v\in C^{\alpha
+\beta }(H)$ and the statement follows.

\subsection{Kolmogorov Equation with Variable Coefficients}

In this section, an estimate is derived to show that $Bu$ is a lower order
operator, which is essential in deriving Schauder estimates in the case of
variable coefficients. To prove Theorem \ref{thm:StoCP}, in a standard way
we use partition of unity and the estimates for constant coefficients, which
allow to obtain apriori estimates. Then the continuation by parameter method
is applied to transfer from constant to variable coefficients.

\subsubsection{Estimates of $Bf$, $f\in C^{\protect\alpha +\protect\beta}$}

\begin{proposition}
\label{b1}Let $\beta \in (0,3)$, $0<\beta \leq \mu <\alpha +\beta $, and 
\begin{equation*}
\int_{|y|\leq 1}|y|^{\alpha }\pi (dy)+\int_{|y|>1}|y|^{\mu }\pi (dy)<\infty .
\end{equation*}%
Assume $a\in C^{\beta }(\mathbf{R}^{d})$,$G^{ij}\in \tilde{C}^{\frac{\beta }{%
\mu \wedge 1}}(\mathbf{R}^{d})$. Then for each $\varepsilon >0$, there
exists a constant $C_{\varepsilon }$ such that%
\begin{equation*}
|Bf|_{\beta }\leq \varepsilon |f|_{\alpha +\beta }+C_{\varepsilon
}|f|_{0},f\in C^{\alpha +\beta }(\mathbf{R}^{d}).
\end{equation*}
\end{proposition}

\begin{proof}
Since the estimates involving the term with $a(x)$ are obvious, in the
following estimates, assume $a=0$.

\textit{Case I}: $\beta \in (0,1]$. Split for $\delta \in (0,1),$%
\begin{eqnarray*}
B_{z}f(x) &=&\int \big[f(x+G(z)y)-f(x)-\mathbf{1}_{\{\alpha \in
(1,2]\}}(\nabla f(x),G(z)y)\chi _{\{|y|\leq 1\}}\big]\pi (dy) \\
&=&\int_{|y|\leq \delta }...+\int_{|y|>\delta
}...=B_{z}^{1}f(x)+B_{z}^{2}f(x)
\end{eqnarray*}%
and $B_{z}^{2}f(x)=B_{z}^{21}f(x)+B_{z}^{22}f(x)$ with 
\begin{eqnarray*}
B_{z}^{21}f(x) &=&f(x)\int_{|y|>\delta }\pi (dy)+\mathbf{1}_{\{\alpha \in
(1,2]\}}\big(\nabla f(x),\int_{\delta <|y|\leq 1}G(z)y\pi (dy)\big), \\
B_{z}^{22}f(x) &=&\int_{|y|>\delta }f(x+G(z)y)\pi (dy).
\end{eqnarray*}%
It follows by the assumptions that there exists $\beta ^{\prime }$ such that 
$\mu <\alpha +\beta ^{\prime }<\alpha +\beta $ and 
\begin{eqnarray}
|B_{z}^{21}f\left( \cdot \right) |_{\beta }+|B_{z}^{22}f\left( \cdot \right)
|_{\beta } &\leq &C[|f|_{\beta }+\mathbf{1}_{\{\alpha \in (1,2]\}}|\nabla
f|_{\beta }],z\in \mathbf{R}^{d},  \label{fo001} \\
|B_{z+h}^{21}f\left( x\right) -B_{z}^{21}f(x)| &\leq &C\mathbf{1}_{\{\alpha
\in (1,2]\}}|\nabla f|_{0}||G||_{\beta }|h|^{\beta },  \notag \\
|B_{z+h}^{22}f(x)-B_{z}^{22}f(x)| &\leq &C\int_{|y|\geq \delta }|y|^{\mu
}\pi (dy)|f|_{\alpha +\beta ^{\prime }}|G|_{\frac{\beta }{\mu \wedge 1}%
},x\in R^{d}.  \notag
\end{eqnarray}

Consider different scenarios on values of $\alpha$ to show that 
\begin{equation}
|B_{z}^{1}f(\cdot )|_{\beta }\leq C|f|_{\alpha +\beta }\int_{|y|\leq \delta
}|y|^{\alpha }d\pi ,z\in \mathbf{R}^{d}.  \label{fo1}
\end{equation}

For $\alpha \in (0,1]$, by Lemma \ref{r1}$,$%
\begin{eqnarray*}
B_{z}^{1}f(x) &=&\int_{|y|\leq \delta }\int \partial ^{\alpha
}f(x-z)k^{(\alpha )}(C(z)y,z)dz\pi (dy)\text{ if }\alpha <1, \\
B_{z}^{1}f(x) &=&\int_{|y|\leq \delta }\int_{0}^{1}(\nabla
f(x+sC(z)y),y)ds\pi (dy)\text{ if }\alpha =1.
\end{eqnarray*}%
Hence, (\ref{fo1}) follows.

For $\alpha =2$, (\ref{fo1}) follows since 
\begin{equation}
B_{z}^{1}f(x)=\int_{|y|\leq \delta }\Big[\int_{0}^{1}\Big( %
D^{2}f(x+sC(z)y)C(z)y,C(z)y\Big) (1-s)ds\Big]\text{ if }\alpha =2.  \notag
\end{equation}

For $\alpha \in (1,2)$, (\ref{fo1}) follows since by Lemma \ref{r1} 
\begin{eqnarray*}
B_{z}^{1}f(x) &=&\int_{|y|\leq \delta }\Big[\int_{0}^{1}\Big(\nabla
f(x+sC(z)y)-\nabla f(x),C(z)y\Big)ds\Big]d\pi \\
&=&\int_{|y|\leq \delta }\Big[\int_{0}^{1}\Big(\int \partial ^{\alpha
-1}\nabla f(x-t)k^{(\alpha -1)}(sC(z)y,t)dt,C(z)y\Big)ds\Big]d\pi .
\end{eqnarray*}

Similarly, to estimate $|B_{\cdot }^{1}f(x)|_{\beta }$, consider different
scenarios on values of $\alpha$.

For $\alpha \in (0,1)$,%
\begin{equation}
|B_{\cdot }^{1}f(x)|_{\beta }\leq |f|_{(\alpha +\beta )}\int_{|y|\leq \delta
}|y|^{\alpha }\pi (dy)||G||_{\frac{\beta }{\mu \wedge 1}},x\in \mathbf{R}%
^{d}.  \label{fo03}
\end{equation}

For $\alpha \in \lbrack 1,2]$, let $\beta \leq \mu <\alpha +\beta ^{\prime
}<\alpha +\beta$.

If $\alpha \in (1,2]$, for $|y|\leq 1,z,\bar{z}\in \mathbf{R}^{d},$%
\begin{eqnarray*}
&&|[f(x+G(z)y)-f(x)-(\nabla f(x),G(z)y)] \\
&&-\left[ f(x+G(\bar{z})y)-f(x)-(\nabla f(x),G(\bar{z})y)\right] | \\
&=&|[f(x+G(z)y)-f(x+G(\bar{z})y)-(\nabla f(x),(G(z)-G(\bar{z}))y)]| \\
&\leq &\int_{0}^{1}|\left( \nabla f(x+(1-s)G(\bar{z})y+sG(z)y)-\nabla
f(x),G(z)y-G(\bar{z})y\right) ds| \\
&\leq &|f|_{\alpha }\left( |G|_{0}^{\alpha -1}|y|^{\alpha -1}+|G(\bar{z}%
)-G(z)|^{\alpha -1}|y|^{\alpha -1}\right) |G(z)-G(\bar{z})||y| \\
&\leq &C|f|_{\alpha }|G|_{0}^{\alpha -1}|G(\bar{z})-G(z)||y|^{\alpha }
\end{eqnarray*}%
and if $\alpha =1,$%
\begin{eqnarray*}
&&|[f(x+G(z)y)-f(x)]-\left[ f(x+G(\bar{z})y)-f(x)\right] | \\
&\leq &\int_{0}^{1}|\left( \nabla f(x+(1-s)G(\bar{z})y+sG(z)y),G(z)y-G(\bar{z%
})y\right) ds| \\
&\leq &|\nabla f|_{0}|G(z)-G(\bar{z})||y|.
\end{eqnarray*}%
It then follows: 
\begin{equation}
|B_{z+h}^{1}f(x)-B_{z}^{1}f(x)|\leq C|h|^{\beta }|f|_{\alpha +\beta ^{\prime
}}||G||_{\beta }^{\alpha },x\in R^{d}.  \label{fo3}
\end{equation}

By (\ref{fo001})-(\ref{fo3}), for each $\varepsilon >0$ there exists a
constant $C_{\varepsilon }$ such that%
\begin{equation}
|Bf|_{\beta }\leq \varepsilon |f|_{\alpha +\beta }+C_{\varepsilon }|f|_{0}.
\label{fos}
\end{equation}

\textit{Case II}: $\beta \in (1,2],$ $\beta \leq \mu <\alpha +\beta .$ Note
that%
\begin{equation*}
\partial _{j}(Bf(x))=\big(\frac{\partial }{\partial _{z_{j}}}B_{z}f(x)\big)%
|_{z=x}+B_{z}f_{x_{j}}|_{z=x}=\big(\frac{\partial }{\partial _{z_{j}}}%
B_{z}f(x)\big)|_{z=x}+Bf_{x_{j}}.
\end{equation*}%
For the second term, apply estimate (\ref{fos}) of \textit{Case I}: $%
f_{x_{j}}\in C^{\alpha +\beta -1},$ the tail moment is 1 and $\beta -1$ $%
\leq \mu -1<\alpha +\beta -1$. Also note that $|G|_{\frac{\beta -1}{(\mu
-1)\wedge 1}}\leq C|G|_{\beta }<\infty $. Hence,%
\begin{equation*}
|Bf_{x_{j}}|_{\beta -1}\leq \varepsilon |f_{x_{j}}|_{\alpha +\beta
-1}+C_{\varepsilon }|f_{x_{j}}|_{0}.
\end{equation*}%
Only the first term needs to be estimated: 
\begin{eqnarray}
B_{z}^{j}f(x) &=&\frac{\partial }{\partial z_{j}}B_{z}f(x)  \label{fo7} \\
&=&\int \left[ \nabla f(x+G(z)y)G_{z_{j}}(z)y-\mathbf{1}_{\{\alpha \in
(1,2]\}}\nabla f(x)G_{z_{j}}(z)y\chi _{\left\{ |y|\leq 1\right\} }\right]
d\pi .  \notag
\end{eqnarray}

Let $B^{j}f(x)=B_{z}^{j}f(x)|_{z=x}$, $x\in R^{d}$. Consider different
scenarios on values of $\alpha$ to show that for each $\varepsilon >0$ there
exists a constant $C_{\varepsilon }$ such that%
\begin{equation}
|B^{j}f|_{\beta -1}\leq \varepsilon |f|_{\alpha +\beta }+C_{\varepsilon
}|f|_{0},f\in C^{\alpha +\beta }(\mathbf{R}^{d}).  \label{fo8}
\end{equation}

For $\alpha \in (0,1]$, since 
\begin{equation*}
B_{z}^{j}f(x)=\int \nabla f(x+G(z)y)G_{z_{j}}(z)yd\pi ,
\end{equation*}%
then for $\mu <\alpha +\beta ^{\prime }<\alpha +\beta $, 
\begin{eqnarray*}
|B_{z}^{j}f(\cdot )|_{\beta -1} &\leq &C|\nabla f|_{\beta
-1}|G_{z_{j}}|_{0},z\in \mathbf{R}^{d}, \\
|B_{z+h}^{j}f(x)-B_{z}^{j}f(x)| &\leq &C|h|^{\beta -1}[|\nabla
f|_{0}||G_{z_{j}}||_{\beta -1}+|f|_{\alpha +\beta ^{\prime
}}|G_{z_{j}}|_{0}||G||_{\frac{\beta }{\mu \wedge 1}}],x\in \mathbf{R}^{d}.
\end{eqnarray*}

For $\alpha \in (1,2]$, split%
\begin{equation*}
B_{z}^{j}f(x)=\int_{|y|\leq
1}...+\int_{|y|>1}...=B_{z}^{j,1}f(x)+B_{z}^{j,2}f(x).
\end{equation*}%
Since by Lemma \ref{r1} 
\begin{equation*}
B_{z}^{j,1}f(x)=\int_{|y|\leq 1}\int \partial ^{\alpha -1}\nabla
f(x-t)k^{(\alpha -1)}(G(z)y,t)G_{z_{j}}(z)ydtd\pi ,
\end{equation*}%
then $\displaystyle|B_{z}^{j,1}f(\cdot )|_{\beta -1}\leq C|f|_{\alpha +\beta
^{\prime }}||G||_{\beta }^{\alpha },z\in \mathbf{R}^{d}$, for some $\beta
^{\prime }\in (0,\beta )$. \newline
For $|y|\leq 1$, $z,\bar{z}\in \mathbf{R}^{d},$ 
\begin{eqnarray*}
&&[\nabla f(x+G(z)y)-\nabla f(x)]G_{z_{j}}(z)y]-[\nabla f(x+G(\bar{z}%
)y)-\nabla f(x)]G_{z_{j}}(\bar{z})y] \\
&\leq &|\nabla f(x+G(z)y)-\nabla f(x+G(\bar{z})y)|~|G_{z_{j}}(z)y| \\
&&+|\nabla f(x+sG(\bar{z})y)-\nabla f(x)||G_{z_{j}}(z)-G_{z_{j}}(\bar{z})||y|
\\
&\leq &|D^{2}f|_{0}|y|^{2}[|G|_{0}|G_{z_{j}}(z)-G_{z_{j}}(\bar{z}%
)|+|G_{z_{j}}|_{0}|G(z)-G(\bar{z})|]
\end{eqnarray*}%
and 
\begin{equation*}
|B_{z}^{j,1}f-B_{\bar{z}}^{j,1}f|\leq C|D^{2}f|_{0}||G||_{\beta }^{2}|z-\bar{%
z}|^{\beta -1}.
\end{equation*}%
Since%
\begin{equation*}
B_{z}^{j,2}f(x)=\int_{|y|>1}\nabla f(x+G(z)y)G_{z_{j}}(z)yd\pi ,
\end{equation*}%
then 
\begin{eqnarray*}
|B_{z}^{j,2}f(\cdot )|_{\beta -1} &\leq &C|\nabla f|_{\beta
-1}|G_{z_{j}}|_{0}, \\
|B_{z+h}^{j,2}f(x)-B_{z}^{j,2}f(x)| &\leq &C|h|^{\beta -1}[|\nabla
f|_{0}||G_{z_{j}}||_{\beta -1}+\int_{|y|>1}|y|^{\mu }d\pi
|D^{2}f|_{0}(1+||G||_{\beta }^{2})],
\end{eqnarray*}%
$z,x,h\in \mathbf{R}^{d}.$

It hence proves that (\ref{fo8}) holds.

\textit{Case III}: $\beta \in (2,3),$ $\beta \leq \mu <\alpha +\beta $. Since%
\begin{equation*}
\partial _{j}(Bf(x))=\big(\frac{\partial }{\partial _{z_{j}}}B_{z}f(x)\big)%
|_{z=x}+B_{z}f_{x_{j}}|_{z=x},
\end{equation*}%
then 
\begin{eqnarray}
\frac{\partial ^{2}}{\partial x_{i}\partial x_{j}}\left( Bf(x)\right)
&=&B_{z}f_{x_{i}x_{j}}(x)|_{z=x}+\frac{\partial }{\partial z_{i}}\left(
B_{z}f_{x_{j}}\right) |_{z=x}  \label{fo00} \\
&&+\frac{\partial }{\partial _{z_{j}}}B_{z}f_{x_{i}}(x)|_{z=x}+\frac{%
\partial ^{2}}{\partial z_{i}\partial z_{j}}B_{z}f(x)|_{z=x}  \notag \\
&=&B\partial _{ij}^{2}f+B^{i}\partial _{j}f+B^{j}\partial _{i}f+B^{ij}f. 
\notag
\end{eqnarray}%
The estimate (\ref{fos}) of \textit{Case I} can be used for the first term ($%
\beta -2\leq \mu -2<\alpha +\beta -2$ with $\beta -2\in (0,1)$). For each $%
\varepsilon ^{\prime }$, there exists a constant $C_{\varepsilon ^{\prime }}$
such that%
\begin{equation*}
|Bf_{x_{i}x_{j}}|_{\beta -2}\leq \varepsilon ^{\prime }|D^{2}f|_{\alpha
+\beta -2}+C_{\varepsilon ^{\prime }}|D^{2}f|_{0}.
\end{equation*}%
For the second and third term in (\ref{fo00}), estimate (\ref{fo8}) of 
\textit{Case II} is applied. Indeed, $f_{x_{j}}\in C^{\alpha +\beta -1}(%
\mathbf{R}^{d})$, $\beta -1\in (1,2)$, $\beta -1\leq \mu -1<\alpha +\beta
-1. $ Hence, for each $\varepsilon ^{\prime }$, there exists a constant $%
C_{\varepsilon ^{\prime }}$ such that%
\begin{equation*}
|B^{i}f_{x_{j}}|_{\beta -2}+|B^{j}f_{x_{i}}|_{\beta -2}\leq \varepsilon
^{\prime }|\nabla f|_{\alpha +\beta -1}+C_{\varepsilon ^{\prime }}|\nabla
f|_{0}.
\end{equation*}

Therefore, only the last term is new. By (\ref{fo7}),%
\begin{eqnarray*}
\frac{\partial ^{2}}{\partial z_{i}\partial z_{j}}B_{z}f(x) &=&\int
(D^{2}f(x+G(z)y)G_{z_{j}}(z)y,G_{z_{i}}(z)y)d\pi \\
&&+\mathbf{1}_{\{\alpha \in (0,1]\}}\int \nabla
f(x+G(z)y)G_{z_{i}z_{j}}(z)yd\pi \\
&&+\mathbf{1}_{\{\alpha \in (1,2]\}}\int \left( \nabla f(x+G(z)y)-\nabla
f(x),G_{z_{i}z_{j}}(z)y\right) d\pi \\
&=&B_{z}^{ij,1}f(x)+B_{z}^{ij,2}f(x)+B_{z}^{ij,3}f(x),
\end{eqnarray*}%
and for $\alpha \in (1,2]$,%
\begin{equation*}
B_{z}^{ij,3}f(x)=\int
\int_{0}^{1}(D^{2}f(x+sG(z)y)G(z)y,G_{z_{i}z_{j}}(z)y)dsd\pi .
\end{equation*}%
It then follows that for $z\in \mathbf{R}^{d}$, 
\begin{eqnarray*}
|B_{z}^{ij,1}f(\cdot )|_{\beta -2} &\leq &|\partial ^{2}f|_{\beta -2}|\nabla
G|_{0}^{2},|B_{z}^{ij,2}f(\cdot )|_{\beta -2}\leq |\nabla f|_{\beta
-2}|\nabla G|_{0}|\partial ^{2}G|_{0}, \\
|B_{z}^{ij,3}f(\cdot )|_{\beta -2} &\leq &|D^{2}f|_{\beta
-2}|D^{2}G|_{0}|G|_{0}.
\end{eqnarray*}%
Let $\beta \leq \mu <\alpha +\beta ^{\prime }<\alpha +\beta $. Then for $%
x,z,h\in \mathbf{R}^{d}$, 
\begin{eqnarray*}
|B_{z+h}^{ij,1}f(x)-B_{z}^{ij,1}f(x)|_{\beta -2} &\leq &|h|^{\beta
-2}(|D^{2}f|_{0}||\nabla G||_{\beta -2}^{2}+|D^{2}f|_{\alpha +\beta ^{\prime
}}|\nabla G|_{0}^{3}\int |y|^{\mu }d\pi ), \\
|B_{z+h}^{ij,2}f(x)-B_{z}^{ij,2}f(x)|_{\beta -2} &\leq &C|h|^{\beta
-2}(|D^{2}f|_{0}||G||_{\beta }^{2}+|\nabla f|_{0}||G||_{\beta }), \\
|B_{z+h}^{ij,3}f(x)-B_{z}^{ij,3}f(x)|_{\beta -2} &\leq &|h|^{\beta
-2}(|D^{3}f|_{0}||G||_{\beta }^{3}+|D^{2}f|_{0}||G||_{\beta }^{2}).
\end{eqnarray*}

Since%
\begin{eqnarray*}
&&B_{x+h}f(x+h)-B_{x}f(x) \\
&=&B_{x+h}f(x+h)-B_{x+h}f(x) \\
&&+B_{x+h}f(x)-B_{x}f(x)
\end{eqnarray*}%
and%
\begin{eqnarray*}
&&B_{x+h}f(x+h)-2B_{x}f(x)+B_{x-h}f(x-h) \\
&=&B_{x+h}f(x+h)-2B_{x+h}f(x)+B_{x+h}f(x-h) \\
&&+2[B_{x+h}f(x)-B_{x}f(x)+ \\
&&+[B_{x-h}f(x-h)-B_{x+h}f(x-h)],
\end{eqnarray*}%
the statement follows.
\end{proof}

\subsubsection{Proof of Theorem \protect\ref{thm:StoCP}}


The proof follows that of Theorem 5 in~\cite{MiP09}, with some simple
changes.

It is well known that for an arbitrary but fixed $\delta >0$, there exist a
family of cubes $D_{k}\subseteq \tilde{D}_{k}\subseteq \mathbf{R}^{d}$ and a
family of deterministic functions $\eta _{k}\in C_{0}^{\infty }(\mathbf{R}%
^{d})$ with the following properties:

\begin{enumerate}
\item For all $k\geq 1,D_{k}$ and $\tilde{D}_{k}$ have a common center $%
x_{k},$ diam $D_{k}\leq \delta$, dist$(D_{k},\mathbf{R}^{d}\backslash \tilde{%
D}_{k})\leq C\delta $ for a constant $C=C(d)>0$, $\cup _{k}D_{k}=\mathbf{R}%
^{d},$ and $1\leq \sum_{k}\mathbf{1}_{\tilde{D}_{k}}\leq 2^{d}.$

\item For all $k$, $0\leq \eta _{k}\leq 1,\eta _{k}=1$ in $D_{k},\eta _{k}=0$
outside of $\tilde{D}_{k}$ and for all multiindices $\gamma$ with $|\gamma
|\leq 3, $%
\begin{equation*}
|\partial ^{\gamma }\eta _{k}|\leq C(d)\delta ^{-|\gamma |}.
\end{equation*}
\end{enumerate}

For $\alpha \in (0,2)$, $\lambda \geq 0$, $k\geq 1,u\in C^{\alpha +\beta
}(H),$ denote%
\begin{eqnarray*}
Au(t,x) &=&A_{x}u(t,x),\quad Bu(t,x)=B_{x}u(t,x),A_{k}u(t,x)=A_{x_{k}}u(t,x),
\\
E_{k}u(t,x) &=&\int [u(t,x+y)-u(t,x)][\eta _{k}(x+y)-\eta _{k}(x)]m(x_{k},y)%
\frac{dy}{|y|^{d+\alpha }}, \\
F_{k}u(t,x) &=&u(t,x)A_{k}\eta _{k}(x).
\end{eqnarray*}

It is readily checked that there is $\beta ^{\prime }\in (0,\beta )$ and a
constant $C$ such that%
\begin{equation*}
\sup_{k}|E_{k}u(t,\cdot )|_{\beta }\leq C|u|_{\alpha +\beta ^{\prime }}
\end{equation*}%
and%
\begin{equation*}
\sup_{k}|F_{k}u(t,\cdot )|_{\beta }\leq C|u|_{\beta }.
\end{equation*}

Elementary calculation shows that for every $u\in C^{\beta }(H),$%
\begin{eqnarray*}
|u|_{0} &\leq &\sup_{k}\sup_{x}|\eta _{k}(x)u(x)|, \\
|u|_{\beta } &\leq &\sup_{k}|\eta _{k}u|_{\beta }+C|u|_{0}\leq
C\sup_{k}|\eta _{k}u|_{\beta }, \\
\sup_{k}|\eta _{k}u|_{\beta } &\leq &|u|_{\beta }+C|u|_{0}\leq C|u|_{\beta },
\end{eqnarray*}%
with $C=C(\alpha ,\delta ,d)$. In particular, 
\begin{equation}
|u|_{\alpha +\beta }\leq C\sup_{k}|\eta _{k}u|_{\alpha +\beta }.  \label{in1}
\end{equation}

Let $u\in C^{\alpha +\beta }(H)$ be a solution to (\ref{eq1}). Then $\eta
_{k}u$ satisfies the equation%
\begin{eqnarray}
\partial _{t}(\eta _{k}u) &=&A_{k}(\eta _{k}u)-\lambda (\eta _{k}u)+\eta
_{k}(Au-A_{k}u)  \label{eqk} \\
&&+\eta _{k}Bu+\eta _{k}f-F_{k}u-E_{k}u,  \notag
\end{eqnarray}%
and by Proposition \ref{prop1},%
\begin{equation*}
|\eta _{k}u|_{\alpha +\beta }\leq C[|\eta _{k}(Au-A_{k}u)|_{\beta }+|\eta
_{k}Bu|_{\beta }+|\eta _{k}f|_{\beta }+|F_{k}u|_{\beta }+|E_{k}u|_{\beta }.
\end{equation*}%
Hence,%
\begin{equation}
|u|_{\alpha +\beta }\leq C[\sup_{k}|\eta _{k}f|_{\beta }+I],  \label{form00}
\end{equation}%
where%
\begin{eqnarray*}
I &\leq &C_{1}\sup_{k}[|\eta _{k}(Au-A_{k}u)|_{\beta }+|\eta _{k}Bu|_{\beta }
\\
&&+|F_{k}u|_{\beta }+|E_{k}u|_{\beta }]+C_{2}|u|_{0}.
\end{eqnarray*}%
By Lemma \ref{prop2}, there exist $\beta ^{\prime }<\beta $, a constant $C$
not depending on $\delta $ and a constant $C=C(\delta )$ 
\begin{equation*}
|\eta _{k}(Au-A_{k}u)|_{\beta }\leq C[C(\delta )|u|_{\alpha +\beta ^{\prime
}}+\delta ^{\beta }|u|_{\alpha +\beta }].
\end{equation*}%
Therefore, for each $\varepsilon >0$, there is a constant $C=C(\varepsilon ) 
$ such that%
\begin{equation*}
|\eta _{k}(Au-A_{k}u)|_{\beta }\leq \varepsilon |u|_{\alpha +\beta
}+C(\varepsilon )|u|_{0}.
\end{equation*}%
By the estimates of Proposition \ref{b1}, it follows that for each $%
\varepsilon >0$, there exists a constant $C_{\varepsilon }$ such that%
\begin{equation*}
I\leq \varepsilon |u|_{\alpha +\beta }+C_{\varepsilon }|u|_{0}.
\end{equation*}

By (\ref{form00}), 
\begin{equation}
|u|_{\alpha +\beta }\leq C[|f|_{\beta }+|u|_{0}].  \label{form01}
\end{equation}

On the other hand, (\ref{eqk}) holds and by Proposition \ref{prop1},%
\begin{eqnarray*}
|u|_{0} &\leq &\sup_{k}|\eta _{k}u|_{\beta } \\
&\leq &\mu (\lambda )\sup_{k}[|f|_{\beta }+|\eta _{k}(Au-A_{k})|_{\beta
}+|\eta _{k}Bu|_{\beta }+|F_{k}u|_{\beta }+|E_{k}u|_{\beta }],
\end{eqnarray*}%
where $\mu (\lambda )\rightarrow 0$ as $\lambda \rightarrow \infty $. Thus,%
\begin{equation}
|u|_{0}\leq C\mu (\lambda )[|f|_{\beta }+|u|_{\alpha +\beta }].
\label{form02}
\end{equation}%
The inequalities (\ref{form01}) and (\ref{form02}) imply that there exist $%
\lambda _{0}>0$ and a constant $C$ independent of $u$ such that if $\lambda
\geq \lambda _{0}$, 
\begin{equation}
|u|_{\alpha +\beta }\leq C|f|_{\beta }.  \label{form03}
\end{equation}%
In a standard way (see ~\cite{MiP09}), it can be verified that (\ref{form03}%
) holds for all $\lambda \geq 0.$ Again by Proposition \ref{prop1} and (\ref%
{in1}), there exists a constant $C$ such that for all $s\leq t\leq T,$ 
\begin{eqnarray*}
|u(t,\cdot )-u(s,\cdot )|_{\frac{\alpha }{2}+\beta } &\leq &\sup_{k}|\eta
_{k}u(t,\cdot )-\eta _{k}u(s,\cdot )|_{\frac{\alpha }{2}+\beta } \\
&\leq &C(t-s)^{\frac{1}{2}}\left( |f|_{\beta }+|u|_{\alpha +\beta }\right) .
\end{eqnarray*}%
Therefore there exists a constant $C$ such that for all $s\leq t\leq T,$%
\begin{equation*}
|u(t,\cdot )-u(s,\cdot )|_{\frac{\alpha }{2}+\beta }\leq C(t-s)^{\frac{1}{2}%
}|f|_{\beta }.
\end{equation*}

To finish the proof, apply the continuation by parameter argument. Let $\tau
\in \left[ 0,1\right] $, $L_{\tau }u=\tau Lu+\left( 1-\tau \right) \partial
^{\alpha }u$ with $L=A+B$ and introduce the space $\hat{C}^{\alpha +\beta
}\left( H\right) $ of functions $u\in C^{\alpha +\beta }(H)$ such that for
each $\left( t,x\right) $, $u\left( t,x\right) =\int_{t}^{T}F\left(
s,x\right) ds$, where $F\in C^{\beta }\left( H\right) .$ It is a Banach
space with respect to the norm 
\begin{equation*}
\left\vert u\right\vert _{\alpha ,\beta }=\left\vert u\right\vert _{\alpha
+\beta }+\left\vert F\right\vert _{\beta }.
\end{equation*}%
Consider the mappings $T_{\tau }:\hat{C}^{\alpha +\beta }\left( H\right)
\rightarrow C^{\beta }(H)$ defined by $u\left( t,x\right)
=-\int_{t}^{T}F\left( s,x\right) \,ds\longmapsto F+L_{\tau }u$. By Lemma \ref%
{prop2} and Proposition \ref{b1}, for some constant $C$ independent of $\tau
,$ $\left\vert T_{\tau }u\right\vert _{\beta }\leq C\left\vert u\right\vert
_{\alpha ,\beta }$. On the other hand, there exists a constant $C$
independent of $\tau $ such that for all $u\in \hat{C}^{\alpha +\beta
}\left( H\right) $,%
\begin{equation}
\left\vert u\right\vert _{\alpha ,\beta }\leq C\left\vert T_{\tau
}u\right\vert _{\beta }.  \label{cp1}
\end{equation}%
Indeed, 
\begin{equation*}
u\left( t,x\right) =-\int_{t}^{T}F\left( s,x\right) \,ds=\int_{t}^{T}\left(
L_{\tau }u-(F+L_{\tau }u)\right) \,ds.
\end{equation*}%
According to the estimate (\ref{form03}), there exists a constant $C$
independent oh $\tau $ such that 
\begin{equation}
\left\vert u\right\vert _{\alpha +\beta }\leq C\left\vert T_{\tau
}u\right\vert _{\beta }=C\left\vert F+L_{\tau }u\right\vert _{\beta }.
\label{cp2}
\end{equation}%
Hence, by Lemma \ref{prop2}, Proposition \ref{b1} and (\ref{cp2}), 
\begin{eqnarray*}
\left\vert u\right\vert _{\alpha ,\beta } &=&\left\vert u\right\vert
_{\alpha +\beta }+\left\vert F\right\vert _{\beta }\leq \left\vert
u\right\vert _{\alpha +\beta }+\left\vert F+L_{\tau }u\right\vert _{\beta
}+\left\vert L_{\tau }u\right\vert _{\beta } \\
&\leq &C(\left\vert u\right\vert _{\alpha +\beta }+\left\vert F+L_{\tau
}u\right\vert _{\beta })\leq C\left\vert F+L_{\tau }u\right\vert _{\beta
}=C\left\vert T_{\tau }u\right\vert _{\beta },
\end{eqnarray*}%
and (\ref{cp1}) follows. Since $T_{0}$ is an onto map, by Theorem 5.2 in 
\cite{GiT83}, all the $T_{\tau }$ are onto maps and the statement follows.

\subsubsection{Proof of Corollary \protect\ref{lcornew1}}

By Lemma \ref{prop2} and Proposition \ref{b1}, for $g\in C^{\alpha +\beta }(%
\mathbf{R}^{d})$, $|Ag|_{\beta }\leq C|g|_{\alpha +\beta }$ and $|Bg|_{\beta
}\leq C|g|_{\alpha +\beta }$ with a constant $C$ independent of $f$ and $g$.
It then follows from (\ref{eq1}) that there exists a unique solution $\tilde{%
v}\in C^{\alpha +\beta }(H)$ to the Cauchy problem 
\begin{eqnarray}  \label{eqn:cauchy_v}
\big(\partial _{t}+A_{x}+B_{x}\big)\tilde{v}(t,x)
&=&f(t,x)-A_{x}g(x)-B_{x}g(x), \\
\tilde{v}(T,x) &=&0  \notag
\end{eqnarray}%
and $|\tilde{v}|_{\alpha +\beta }\leq C\big(|g|_{\alpha +\beta }+|f|_{\beta }%
\big)$ with $C$ independent of $f$ and $g$. Let $v(t,x)=\tilde{v}(t,x)+g(x)$%
, where $\tilde{v}$ is the solution to problem (\ref{eqn:cauchy_v}). Then $v$
is the unique solution to the Cauchy problem (\ref{maf8}) and $|v|_{\alpha
+\beta }\leq C(|g|_{\alpha +\beta }+|f|_{\beta })$.

\begin{remark}
\label{rlast}If the assumptions of Corollary \ref{lcornew1} hold and $v\in
C^{\alpha +\beta }(H)$ is the solution to $(\ref{maf8})$, then $\partial
_{t}v=f-A_{x}v-B_{x}v$, and by Lemma \ref{prop2} and Proposition \ref{b1}, $%
|\partial _{t}v|_{\beta }\leq C(|g|_{\alpha +\beta }+|f|_{\beta }).$
\end{remark}


\section{One-Step Estimate and Proof of Main Result}


The following Lemma provides a one-step estimate of the conditional
expectation of an increment of the Euler approximation.

\begin{lemma}
\label{lem:expect}Let $\beta \in (0,3)$, $0<\beta \leq \mu <\alpha +\beta $,
and 
\begin{equation*}
\int_{|y|\leq 1}|y|^{\alpha }\pi (dy)+\int_{|y|>1}|y|^{\mu }\pi (dy)<\infty .
\end{equation*}%
Assume $a^{i},b^{ij}\in \tilde{C}^{\beta }(\mathbf{R}^{d}),G^{ij}\in \tilde{C%
}^{\frac{\beta }{\mu \wedge 1}}(\mathbf{R}^{d})$. Then there exists a
constant $C$ such that for all $f\in C^{\beta }(\mathbf{R}^{d}),$%
\begin{equation*}
\big\vert\mathbf{E}\big[f(Y_{s})-f(Y_{\tau _{i_{s}}})|\widetilde{\mathcal{F}}%
_{\tau _{i_{s}}}\big]\big\vert\leq C|f|_{\beta }r(\delta ,\alpha ,\beta
),\forall s\in \lbrack 0,T],
\end{equation*}%
where $i_{s}=i$ if $\tau _{i}\leq s<\tau _{i+1}$ and $r(\delta ,\alpha
,\beta )$ is as defined in Theorem \ref{thm:main}.
\end{lemma}

The proof of Lemma~\ref{lem:expect} is based on applying It\^{o}'s formula
to $f(Y_{s})-f(Y_{\tau _{i_{s}}})$, $f\in C^{\beta }(\mathbf{R}^{d})$. If $%
\beta >\alpha $, by Remark \ref{lrenew1} and It\^{o}'s formula, the
inequality holds. If $\beta \leq \alpha $, $f$ is first smoothed by using $%
w\in C_{0}^{\infty }(\mathbf{R}^{d}),$ a nonnegative smooth function with
support on $\{|x|\leq 1\}$ such that $w(x)=w(|x|)$, $x\in \mathbf{R}^{d},$
and $\int w(x)dx=1$ (see (8.1) in \cite{Fol99})$.$ Note that, due to the
symmetry, 
\begin{equation}
\int_{\mathbf{R}^{d}}x^{i}w(x)dx=0,i=1,\ldots ,d.  \label{ff26}
\end{equation}

For $x\in \mathbf{R}^{d}$ and $\varepsilon \in (0,1)$, define $%
w^{\varepsilon }(x)=\varepsilon ^{-d}w\left( \frac{x}{\varepsilon }\right) $
and the convolution 
\begin{equation}
f^{\varepsilon }(x)=\int f(y)w^{\varepsilon }(x-y)dy=\int
f(x-y)w^{\varepsilon }(y)dy,x\in \mathbf{R}^{d}.  \label{maf7}
\end{equation}


\subsection{Some Auxiliary Estimates}

For the estimates of $A_{z}f^{\varepsilon }$, the following simple integral
estimates are needed. Recall that $m(z,y)$ in the definition of operator $%
A_{z}$ (see (\ref{foo})) is bounded, smooth, and 0-homogeneous and symmetric
in $y$.

\begin{lemma}
\label{rem2} Let $v\in C_{0}^{\infty }(\mathbf{R}^{d})$.

\begin{enumerate}
\item[\textup{(i)}] For $\alpha \in (0,2)$,%
\begin{equation*}
\int_{\mathbf{R}^{d}}\int_{\mathbf{R}_{0}^{d}}|v(y+y^{\prime })-v(y)-\chi
^{(\alpha )}(y^{\prime })(\nabla v(y),y^{\prime })|\frac{dydy^{\prime }}{%
|y^{\prime }|^{d+\alpha }}<\infty ,
\end{equation*}%
where $\chi ^{(\alpha )}(y)=\mathbf{1}_{\{|y|\leq 1\}}\mathbf{1}_{\{\alpha
=1\}}+\mathbf{1}_{\{\alpha \in (1,2)\};}$

\item[\textup{(ii)}] For $\beta \in (0,1]$, $\beta <\alpha ,z\in \mathbf{R}%
^{d},$%
\begin{equation*}
\sup_{z}\int_{\mathbf{R}^{d}}|(A_{z}v)(y)||y|^{\beta }dy<\infty ,
\end{equation*}%
and for $\beta \in (0,1]$, $\beta =\alpha ,z\in \mathbf{R}^{d},k>1,$%
\begin{equation*}
\sup_{z}\int_{\mathbf{R}^{d}}|(A_{z}v)(y)|(|y|^{\alpha }\wedge k)dy\leq
C(1+\ln k).
\end{equation*}

\item[\textup{(iii)}] For $1<\beta <\alpha <2,$%
\begin{equation*}
\int_{\mathbf{R}^{d}}\int_{\mathbf{R}_{0}^{d}}\int_{0}^{1}|v(y+sy^{\prime
})-v(y)|~|y|^{\beta -1}\frac{dsdydy^{\prime }}{|y^{\prime }|^{d+\alpha -1}}%
<\infty ,
\end{equation*}%
and for $1<\beta =\alpha <2,k>1,$%
\begin{eqnarray*}
&&\int_{\mathbf{R}^{d}}\int_{\mathbf{R}_{0}^{d}}\int_{0}^{1}|v(y+sy^{\prime
})-v(y)|~(|y|^{\beta -1}\wedge k)\frac{dsdydy^{\prime }}{|y^{\prime
}|^{d+\alpha -1}} \\
&\leq &C(1+\ln k).
\end{eqnarray*}
\end{enumerate}
\end{lemma}

\begin{proof}
(i) Indeed,%
\begin{eqnarray*}
&&|v(y+y^{\prime })-v(y)-\chi ^{(\alpha )}(y^{\prime })(\nabla
v(y),y^{\prime })| \\
&\leq &\mathbf{1}_{\{|y^{\prime }|\leq 1\}}\big\{\int_{0}^{1}[\max_{i,j}|%
\partial _{ij}^{2}v(y+sy^{\prime })|~|y^{\prime }|^{2}+\mathbf{1}_{\{\alpha
\in (0,1)\}}|\nabla v(y+sy^{\prime })||y^{\prime }|]ds\big\} \\
&&+\mathbf{1}_{\{|y^{\prime }|>1\}}\big\{|v(y+y^{\prime })|+|v(y)|+\mathbf{1}%
_{\{\alpha \in (1,2)\}}|\nabla v(y)|~|y^{\prime }|\big\},y,y^{\prime }\in 
\mathbf{R}^{d}.
\end{eqnarray*}%
The claim follows.

(ii) For $\beta \in (0,1]$, $\beta <\alpha $, $z\in \mathbf{R}^{d}$, 
\begin{eqnarray*}
\int_{\mathbf{R}^{d}}|(A_{z}v)(y)||y|^{\beta }dy &\leq &\int_{\mathbf{R}%
^{d}}\int_{|y^{\prime }|>1}|v(y+y^{\prime })|~|y|^{\beta }\frac{dydy^{\prime
}}{|y^{\prime }|^{d+\alpha }} \\
&&+\int_{\mathbf{R}^{d}}\int_{|y^{\prime }|>1}|v(y)||y|^{\beta }\frac{%
dydy^{\prime }}{|y^{\prime }|^{d+\alpha }} \\
&&+\max_{i,j}\int_{\mathbf{R}^{d}}\int_{|y^{\prime }|\leq
1}\int_{0}^{1}|\partial _{ij}^{2}v(y+sy^{\prime })|~|y^{\prime
}|^{2}|y|^{\beta }\frac{dsdy^{\prime }dy}{|y^{\prime }|^{d+\alpha }}
\end{eqnarray*}%
and%
\begin{eqnarray*}
\int_{\mathbf{R}^{d}}\int_{|y^{\prime }|>1}|v(y+y^{\prime })||y|^{\beta }%
\frac{dydy^{\prime }}{|y^{\prime }|^{d+\alpha }} &\leq &C\Big[\int_{\mathbf{R%
}^{d}}\int_{|y^{\prime }|>1}|v(y+y^{\prime })||y+y^{\prime }|^{\beta }\frac{%
dydy^{\prime }}{|y^{\prime }|^{d+\alpha }} \\
&&+\int_{\mathbf{R}^{d}}\int_{|y^{\prime }|>1}|v(y+y^{\prime })||y^{\prime
}|^{\beta }\frac{dydy^{\prime }}{|y^{\prime }|^{d+\alpha }}\Big].
\end{eqnarray*}%
Let $\beta \in (0,1],\beta =\alpha .$ Assume $v(x)=0$ if $|x|>R$. We have
for $k>1$ with $A=(R+1)^{1/\alpha },$%
\begin{eqnarray*}
&&\int_{\mathbf{R}^{d}}\int_{|y^{\prime }|>1}|v(y+y^{\prime })|~(|y|^{\alpha
}\wedge k)\frac{dydy^{\prime }}{|y^{\prime }|^{d+\alpha }} \\
&\leq &\int_{\mathbf{R}^{d}}\int_{|y^{\prime }|>1}|v(y+y^{\prime
})|~(|y|^{\alpha }\wedge Ak)\frac{dydy^{\prime }}{|y^{\prime }|^{d+\alpha }}
\\
&=&\int_{|y|\leq (R+1)k^{1/\alpha }}\int_{|y^{\prime }|>1}|v(y+y^{\prime
})|~|y|^{\alpha }\frac{dydy^{\prime }}{|y^{\prime }|^{d+\alpha }} \\
&&+k\int_{|y|>(R+1)k^{1/\alpha }}\int_{\left\vert y^{\prime }\right\vert
>1}|v(y+y^{\prime })|\frac{dydy^{\prime }}{|y^{\prime }|^{d+\alpha }} \\
&=&A_{1}+A_{2.}
\end{eqnarray*}%
Then 
\begin{eqnarray*}
|A_{1}| &\leq &\int \int_{1\leq |y^{\prime }|\leq (R+1)(1+k^{1/\alpha
})}|v(y+y^{\prime })|\frac{dy^{\prime }dy}{|y^{\prime }|^{d}} \\
&\leq &C(1+\ln k).
\end{eqnarray*}%
Since for $|y+y^{\prime }|\leq R,|y|>(R+1)k^{1/\alpha },$ we have $%
|y^{\prime }|\geq (R+1)k^{1/\alpha }-R\geq k^{1/\alpha }$, it follows 
\begin{eqnarray*}
|A_{2}| &\leq &k\int \int_{|y^{\prime }|\geq k^{1/\alpha }}|v(y+y^{\prime })|%
\frac{dydy^{\prime }}{|y^{\prime }|^{d+\alpha }} \\
&\leq &Ckk^{-1}=C.
\end{eqnarray*}%
Then%
\begin{eqnarray*}
&&\int \int_{|y^{\prime }|\leq 1}|v(y+y^{\prime })-v(y)-1_{\alpha =1}(\nabla
v(y),y^{\prime })|~(|y|^{\alpha }\wedge k)\frac{dydy^{\prime }}{|y^{\prime
}|^{d+\alpha }} \\
&\leq &\int_{0}^{1}\int \int_{|y^{\prime }|\leq 1}|\nabla v(y+sy^{\prime
})-1_{\alpha =1}\nabla v(y)|~|y|^{\alpha }\frac{dydy^{\prime }}{|y^{\prime
}|^{d+\alpha -1}}<\infty .
\end{eqnarray*}%
Part (ii) follows.

(iii) For $1<\beta <\alpha <2,$ 
\begin{eqnarray*}
&&\int_{\mathbf{R}^{d}}\int_{\mathbf{R}_{0}^{d}}\int_{0}^{1}|v(y+sy^{\prime
})-v(y)|~|y|^{\beta -1}\frac{dydy^{\prime }ds}{|y^{\prime }|^{d+\alpha -1}}
\\
&\leq &\int_{\mathbf{R}^{d}}\int_{|y^{\prime
}|>1}\int_{0}^{1}|v(y+sy^{\prime })||y|^{\beta -1}\frac{dydy^{\prime }ds}{%
|y^{\prime }|^{d+\alpha -1}} \\
&&+\int_{\mathbf{R}^{d}}\int_{|y^{\prime }|>1}\int_{0}^{1}|v(y)||y|^{\beta
-1}\frac{dydy^{\prime }ds}{|y^{\prime }|^{d+\alpha -1}} \\
&&+\int_{\mathbf{R}^{d}}\int_{|y^{\prime }|\leq
1}\int_{0}^{1}\int_{0}^{1}|\nabla v(y+s\tau y^{\prime })||y|^{\beta -1}\frac{%
dsd\tau dydy^{\prime }}{|y^{\prime }|^{d+\alpha -2}}.
\end{eqnarray*}%
Since%
\begin{eqnarray*}
&&\int_{\mathbf{R}^{d}}\int_{|y^{\prime }|>1}\int_{0}^{1}|v(y+sy^{\prime
})||y|^{\beta -1}\frac{dydy^{\prime }ds}{|y^{\prime }|^{d+\alpha -1}} \\
&\leq &C\Big[\int_{\mathbf{R}^{d}}\int_{|y^{\prime
}|>1}\int_{0}^{1}|v(y+sy^{\prime })||y+sy^{\prime }|^{\beta -1}\frac{%
dydy^{\prime }ds}{|y^{\prime }|^{d+\alpha -1}} \\
&&+\int_{\mathbf{R}^{d}}\int_{|y^{\prime }|>1}\int_{0}^{1}|v(y+sy^{\prime
})||y^{\prime }|^{\beta -1}\frac{dydy^{\prime }ds}{|y^{\prime }|^{d+\alpha
-1}}\Big]
\end{eqnarray*}%
and similarly%
\begin{eqnarray*}
&&\int_{\mathbf{R}^{d}}\int_{|y^{\prime }|\leq
1}\int_{0}^{1}\int_{0}^{1}|\nabla v(y+s\tau y^{\prime })||y|^{\beta -1}\frac{%
dsd\tau dydy^{\prime }}{|y^{\prime }|^{d+\alpha -2}} \\
&\leq &C\Big[\int_{\mathbf{R}^{d}}\int_{|y^{\prime }|\leq
1}\int_{0}^{1}\int_{0}^{1}|\nabla v(y+s\tau y^{\prime })||y+s\tau y^{\prime
}|^{\beta -1}\frac{dsd\tau dydy^{\prime }}{|y^{\prime }|^{d+\alpha -2}} \\
&&+\int_{\mathbf{R}^{d}}\int_{|y^{\prime }|\leq
1}\int_{0}^{1}\int_{0}^{1}|\nabla v(y+s\tau y^{\prime })||y^{\prime
}|^{\beta -1}\frac{dsd\tau dydy^{\prime }}{|y^{\prime }|^{d+\alpha -2}}\Big]
\end{eqnarray*}%
are finite, the first estimate in (iii) follows.

If $1<\beta =\alpha <2$, we simply repeat the proof in part (ii). The
statement follows.
\end{proof}

\begin{remark}
The estimate in part (iii) implies that (ii) can be extended to all $0<\beta
\leq \alpha <2.$
\end{remark}

We will use the following modulus of continuity estimate of a function $f\in
C^{1}(\mathbf{R}^{d}).$

\begin{lemma}
\label{lemn1}(cf. Lemma 5.6 in \cite{Caf}, Lemma 2.2 in \cite{kimd}) Let $%
f\in C^{1}(\mathbf{R}^{d})$ and $[f]_{1}\leq K$. Then there is a constant $C$
such that for all $x,h\in \mathbf{R}^{d},h\neq 0,$%
\begin{equation*}
|f(x+h)-f(x)|\leq C|h|(1+\left\vert \ln |h|\right\vert )|f|_{1}.
\end{equation*}
\end{lemma}

\begin{proof}
We follow the steps of Lemma 5.6 in \cite{Caf}. Fix $x,h\in \mathbf{R}%
^{d},h\neq 0$ such that $0<|h|<1/2$, we find a positive integer $k$ so that%
\begin{equation*}
2^{-k-1}\leq |h|<2^{-k}
\end{equation*}%
or $\frac{1}{2}2^{-k}\leq |h|<2^{-k}$. Set $\tau _{0}=2^{k}h$ (\textbf{note: 
}$\frac{1}{2}\leq \tau _{0}<1,2^{-k}\leq 2|h|),\ln |h|<-k\ln 2$ or $k<\frac{%
-\ln |h|}{\ln 2}$). Define for $\tau \in \mathbf{R}^{d},$%
\begin{equation*}
v(\tau )=f(x+\tau )-f(x).
\end{equation*}%
Note that%
\begin{eqnarray*}
|v(\tau )-2v(\tau /2)| &=&|f(x+\tau )-2f(x+\tau /2)+f(x)| \\
&\leq &[f]_{1}|\tau |/2.
\end{eqnarray*}%
Thus, 
\begin{equation*}
|2^{j-1}v(\tau _{0}/2^{j-1})-2^{j}v(\tau _{0}/2^{j})|\leq \lbrack
f]_{1}2^{j-1}|\tau _{0}|/2^{j}=2^{-1}[f]_{1}|\tau _{0}|,
\end{equation*}%
which implies ($v(\tau _{0})-2^{k}v(h)=\sum_{j=1}^{k}\left( 2^{j-1}v(\tau
_{0}/2^{j-1}\right) -\left( 2^{j}v(\tau _{0}/2^{j}\right) )$ 
\begin{eqnarray*}
|v(\tau _{0})-2^{k}v(h)| &\leq &\sum_{j=1}^{k}|2^{j-1}v(\tau
_{0}/2^{j-1})-2^{j}v(\tau _{0}/2^{j})| \\
&\leq &k2^{-1}|\tau _{0}|[f]_{1}.
\end{eqnarray*}%
Since $|v(\tau _{0})|\leq 2|f|_{0}$ or $|v(\tau _{0})|\leq \lbrack f]_{1/2})$%
, we derive%
\begin{eqnarray*}
|v(h)| &\leq &2^{-k}|2^{k}v(h)-v(\tau _{0})|+2^{-k}|v(\tau _{0})| \\
&\leq &[f]_{1}k2^{-1-k}|\tau _{0}|+2\cdot 2^{-k}|f|_{0} \\
&\leq &[f]_{1}k|h|+4|h||f|_{0} \\
&\leq &C|f|_{1}|h|(1-\ln |h|).
\end{eqnarray*}%
The statement follows.
\end{proof}

Now we prove some estimates for $Af^{\varepsilon }$ and $Bf^{\varepsilon }$.

\begin{lemma}
\label{lnew2}Let $\varepsilon \in (0,1).$

(i) Let $\alpha \in (0,2)$. Then there exists a constant $C$ such that for
all $z,x\in \mathbf{R}^{d},$ 
\begin{equation}
|A_{z}f^{\varepsilon }(x)|\leq C\kappa (\varepsilon ,\alpha ,\beta
)|f|_{\beta },  \label{ff29}
\end{equation}%
where $\kappa (\varepsilon ,\alpha ,\beta )=\varepsilon ^{-\alpha +\beta }$
if $\beta <\alpha $ and $\kappa (\varepsilon ,\alpha ,\beta )=1-\ln
\varepsilon $ if $\beta =\alpha ;$ in particular, for all $f\in C^{\beta }(%
\mathbf{R}^{d}),z,x\in \mathbf{R}^{d}$, 
\begin{equation}
|\partial ^{\alpha }f^{\varepsilon }(x)|\leq C\kappa (\varepsilon ,\alpha
,\beta )|f|_{\beta }.  \label{maf5}
\end{equation}

(ii)For each $\beta \in (0,2]$ there exists a constant $C$ such that for all 
$f\in C^{\beta }(\mathbf{R}^{d}),x\in \mathbf{R}^{d}$, 
\begin{equation*}
|f^{\varepsilon }(x)-f(x)|\leq C\gamma (\varepsilon ,\beta )|f|_{\beta },
\end{equation*}%
where $\gamma (\varepsilon ,\beta )=\varepsilon ^{\beta }$ if $\beta <2$ and 
$\gamma (\varepsilon ,2)=\varepsilon ^{2}(1-\ln \varepsilon ).$

(iii) Let $\beta \in (0,2]$. For $k,l=1,\ldots ,d,x\in \mathbf{R}^{d},$%
\begin{eqnarray}
|\partial _{k}f^{\varepsilon }(x)| &\leq &C\kappa (\varepsilon ,1,\beta
)|f|_{\beta },\mbox{ if }\beta \leq 1,  \label{ff30} \\
|\partial _{kl}^{2}f^{\varepsilon }(x)| &\leq &C\kappa (\varepsilon ,2,\beta
)|f|_{\beta },\mbox{ if }\beta \leq 2,  \notag \\
|f^{\varepsilon }|_{1} &\leq &C|f|_{1},  \notag
\end{eqnarray}%
and%
\begin{eqnarray}
|f^{\varepsilon }|_{\alpha } &\leq &C\varepsilon ^{-\alpha +\beta
}|f|_{\beta },\mbox{ if }\beta \in (0,1],\alpha \in \lbrack 1,2),
\label{maf5'} \\
|\partial ^{\alpha -1}\nabla f^{\varepsilon }(x)| &\leq &C\kappa
(\varepsilon ,\alpha ,\beta )|f|_{\beta },\mbox{ if }\beta \in (0,\alpha
],\alpha \in (1,2),\beta \neq \alpha -1.  \label{maf6}
\end{eqnarray}
\end{lemma}

\begin{proof}
(i) For $z,x\in \mathbf{R}^{d},$ by changing the variable of integration
with $\bar{y}=\frac{y}{\varepsilon }$ and using (\ref{ff27}) for $\alpha =1$%
, 
\begin{eqnarray}
A_{z}w^{\varepsilon }(x) &=&\mathbf{1}_{\{\alpha =1\}}(a(z),\nabla
w^{\varepsilon }(x))  \notag \\
&&+\int [w^{\varepsilon }(x+y)-w^{\varepsilon }(x)-\bar{\chi}_{\alpha
}(y)(\nabla w^{\varepsilon }(x),y)]m(z,y)\frac{dy}{|y|^{d+\alpha }}  \notag
\\
&=&\varepsilon ^{-\alpha }\varepsilon ^{-d}(A_{z}w)(\frac{x}{\varepsilon }),
\label{and}
\end{eqnarray}%
where $\bar{\chi}_{\alpha }(y)=\mathbf{1}_{\{|y|\leq 1\}}\mathbf{1}%
_{\{\alpha =1\}}+\mathbf{1}_{\{\alpha \in (1,2)\}},y\in \mathbf{R}^{d}.$ It
follows from Lemma \ref{rem2}(i), the Fubini theorem, and (\ref{and}),
changing the variable of integration with $\bar{y}=\frac{y}{\varepsilon }$
as well, that%
\begin{eqnarray*}
A_{z}f^{\varepsilon }(x) &=&\int_{\mathbf{R}^{d}}\varepsilon ^{-\alpha
}\varepsilon ^{-d}(A_{z}w)(\frac{x-y}{\varepsilon })f(y)dy \\
&=&\int \varepsilon ^{-\alpha }\varepsilon ^{-d}(A_{z}w)(\frac{y}{%
\varepsilon })f(x-y)dy \\
&=&\int \varepsilon ^{-\alpha }(A_{z}w)(y)f(x-\varepsilon y)dy,x,z\in 
\mathbf{R}^{d}.
\end{eqnarray*}

By Lemma \ref{rem2}(i) and the Fubini theorem,%
\begin{equation*}
\int_{\mathbf{R}^{d}}A_{z}w(y)dy=0.
\end{equation*}%
Also, it is easy to see that $A_{z}w(y)=A_{z}w(-y),y\in \mathbf{R}^{d}$.
Hence, if $\beta \in (0,1]$, $\beta \leq \alpha $,%
\begin{eqnarray*}
A_{z}f^{\varepsilon }(x) &=&\int \varepsilon ^{-\alpha
}(A_{z}w)(y)f(x-\varepsilon y)dy \\
&=&\frac{1}{2}\int \varepsilon ^{-\alpha }(A_{z}w)(y)[f(x-\varepsilon
y)+f(x+\varepsilon y)-2f(x)]dy
\end{eqnarray*}%
and 
\begin{equation*}
|A_{z}f^{\varepsilon }(x)|\leq C\varepsilon ^{-\alpha +\beta }|f|_{\beta
}\int_{\mathbf{R}^{d}}|(A_{z}w)(y)|~(|y|^{\beta }\wedge \varepsilon ^{-\beta
})dy.
\end{equation*}%
So, by Lemma \ref{rem2}, (\ref{ff29}) holds for $\beta \leq \alpha ,\beta
\in (0,1].$

Assume $1<\beta \leq \alpha <2$. By Theorem 2.27 in \cite{Fol99},
differentiation and integration can be can switched:%
\begin{eqnarray*}
A_{z}w(y) &=&\int [w(y+y^{\prime })-w(y)-(\nabla w(y),y^{\prime
})]m(z,y^{\prime })\frac{dy^{\prime }}{|y^{\prime }|^{d+\alpha }} \\
&=&\int \int_{0}^{1}\big(\nabla _{y}w(y+sy^{\prime })-\nabla
_{y}w(y),y^{\prime }\big)dsm(z,y^{\prime })\frac{dy^{\prime }}{|y^{\prime
}|^{d+\alpha }} \\
&=&\sum_{i=1}^{d}\frac{\partial }{\partial y_{i}}\int
\int_{0}^{1}[w(y+sy^{\prime })-w(y)]y_{i}^{\prime }dsm(z,y^{\prime })\frac{%
dy^{\prime }}{|y^{\prime }|^{d+\alpha }}.
\end{eqnarray*}%
By integrating by parts, 
\begin{eqnarray}
A_{z}f^{\varepsilon }(x) &=&\int \varepsilon ^{-\alpha
}A_{z}w(y)f(x-\varepsilon y)dy  \notag \\
&=&\varepsilon ^{-\alpha +1}\int \int_{0}^{1}[w(y+sy^{\prime })-w(y)]
\label{ff28} \\
&&\times (\nabla f(x-\varepsilon y),y^{\prime })m(z,y^{\prime })\frac{%
dsdydy^{\prime }}{|y^{\prime }|^{d+\alpha }},x\in \mathbf{R}^{d}.  \notag
\end{eqnarray}%
Since 
\begin{equation*}
\int_{\mathbf{R}_{0}^{d}}\int_{\mathbf{R}^{d}}\int_{0}^{1}|w(y+sy^{\prime
})-w(y)|~|y^{\prime }|\frac{dsdydy^{\prime }}{|y^{\prime }|^{d+\alpha }}%
<\infty ,
\end{equation*}%
the Fubini theorem applies, $\int [w(y+sy^{\prime })-w(y)]dy=0$ and we can
rewrite (\ref{ff28}) as%
\begin{eqnarray*}
A_{z}f^{\varepsilon }(x) &=&\varepsilon ^{-\alpha +1}\int_{\mathbf{R}%
^{d}}\int_{\mathbf{R}_{0}^{d}}\int_{0}^{1}[w(y+sy^{\prime })-w(y)] \\
&&\times (\nabla f(x-\varepsilon y)-\nabla f(x),y^{\prime })m(z,y^{\prime })%
\frac{dsdydy^{\prime }}{|y^{\prime }|^{d+\alpha }},x,z\in \mathbf{R}^{d}.
\end{eqnarray*}%
Hence, 
\begin{eqnarray*}
|A_{z}f^{\varepsilon }(x)| &\leq &C\varepsilon ^{-\alpha +1}\varepsilon
^{\beta -1}|\nabla f|_{\beta -1}\int \int \int_{0}^{1}|w(y+sy^{\prime
})-w(y)| \\
&&\times (|y|^{\beta -1}\wedge \varepsilon ^{-(\beta -1)})\frac{dydy^{\prime
}}{|y^{\prime }|^{d+\alpha -1}},
\end{eqnarray*}%
and by Lemma \ref{rem2}(iii), (\ref{ff29}) is proved for $1<\beta \leq
\alpha <2$. By taking $m=1,$ (\ref{maf5}) follows.

(ii) For $\beta \in (1,2)$, by (\ref{ff26}), 
\begin{eqnarray*}
f^{\varepsilon }(x)-f(x) &=&\int [f(x-y)-f(x)]w^{\varepsilon }(y)dy \\
&=&\int [f(x+y)-f(x)-(\nabla f(x),y)]w^{\varepsilon }(y)dy \\
&=&\int \int_{0}^{1}(\nabla f(x+sy)-\nabla f(x),y)dsw^{\varepsilon }(y)dy
\end{eqnarray*}

and%
\begin{equation*}
|f^{\varepsilon }(x)-f(x)|\leq C|\nabla f|_{\beta -1}\int |y|^{1+(\beta
-1)}w^{\varepsilon }(y)dy\leq C|f|_{\beta }\varepsilon ^{\beta }.
\end{equation*}%
If $\beta =2$, then by Lemma \ref{lemn1},%
\begin{eqnarray*}
|f^{\varepsilon }(x)-f(x)| &\leq &C|\nabla f|_{1}\int |y|^{2}(1+\left\vert
\ln |y|\right\vert )w^{\varepsilon }(y)dy \\
&\leq &C|\nabla f|_{1}\varepsilon ^{2}(1-\ln \varepsilon ).
\end{eqnarray*}

For $\beta \in (0,1],$%
\begin{eqnarray*}
f^{\varepsilon }(x)-f(x) &=&\int [f(x-y)-f(x)]w^{\varepsilon }(y)dy \\
&=&\int [f(x+y)-f(x)]w^{\varepsilon }(y)dy
\end{eqnarray*}%
and%
\begin{equation*}
f^{\varepsilon }(x)-f(x)=\frac{1}{2}\int [f(x+y)+f(x-y)-2f(x)]w^{\varepsilon
}(y)dy.
\end{equation*}%
Hence, for $\beta \in (0,1],$%
\begin{equation*}
|f^{\varepsilon }(x)-f(x)|\leq C|f|_{\beta }\varepsilon ^{\beta }.
\end{equation*}

(iii) If $\beta <1$, by changing the variable of integration,%
\begin{eqnarray*}
\partial _{k}f^{\varepsilon }(x) &=&\varepsilon ^{-1}\int_{\mathbf{R}%
^{d}}\varepsilon ^{-d}\partial _{k}w(\frac{x-y}{\varepsilon })f(y)dy \\
&=&\varepsilon ^{-1}\int_{\mathbf{R}^{d}}\varepsilon ^{-d}\partial _{k}w(%
\frac{y}{\varepsilon })f(x-y)dy \\
&=&\varepsilon ^{-1}\int_{\mathbf{R}^{d}}\partial _{k}w(y)[f(x-\varepsilon
y)-f(x)]dy,
\end{eqnarray*}%
and the first inequality follows by Lemma \ref{lemn1}. Since%
\begin{eqnarray*}
&&f^{\varepsilon }(x+h)+f^{\varepsilon }(x-h)-2f^{\varepsilon }(x) \\
&=&\frac{1}{2}\int w_{\varepsilon }(y)[f(x-y+h)+f(x-y-h)-2f(x-y)]dy,
\end{eqnarray*}%
we have $|f^{\varepsilon }|_{1}\leq |f|_{1}$. Also, since $\partial
_{kl}^{2}w(y)=\partial _{kl}^{2}w(-y),k,l=1,\ldots ,d,y\in \mathbf{R}^{d},$ 
\begin{eqnarray*}
\partial _{kl}^{2}f^{\varepsilon }(x) &=&\varepsilon ^{-2}\int_{\mathbf{R}%
^{d}}\varepsilon ^{-d}\partial _{kl}^{2}w(\frac{x-y}{\varepsilon })f(y)dy \\
&=&\varepsilon ^{-2}\int_{\mathbf{R}^{d}}\varepsilon ^{-d}\partial
_{kl}^{2}w(\frac{y}{\varepsilon })f(x-y)dy \\
&=&\varepsilon ^{-2}\int_{\mathbf{R}^{d}}\partial
_{kl}^{2}w(y)[f(x-\varepsilon y)-f(x)]dy \\
&=&\frac{1}{2}\varepsilon ^{-2}\int_{\mathbf{R}^{d}}\partial
_{kl}^{2}w(y)[f(x+\varepsilon y)+f(x-\varepsilon y)-2f(x)]dy.
\end{eqnarray*}

\ \ \ \ \ Thus, for all $x\in \mathbf{R}^{d}$, 
\begin{equation*}
\ |\partial _{kl}^{2}f^{\varepsilon }(x)|\leq C\varepsilon ^{-2+\beta
}|f|_{\beta }\text{ if }\beta \in (0,1].
\end{equation*}%
Similarly, if $1<\beta \leq 2,$%
\begin{eqnarray*}
\partial _{k}f^{\varepsilon }(x) &=&\int \varepsilon ^{-d}w(\frac{y}{%
\varepsilon })\partial _{k}f(x-y)dy \\
&=&\int \varepsilon ^{-d}w(\frac{x-y}{\varepsilon })\partial _{k}f(y)dy
\end{eqnarray*}%
and%
\begin{eqnarray*}
\partial _{kl}^{2}f^{\varepsilon }(x) &=&\varepsilon ^{-1}\int \varepsilon
^{-d}\partial _{l}w(\frac{y}{\varepsilon })\partial _{k}f(x-y)dy \\
&=&\varepsilon ^{-1}\int \partial _{l}w(y)[\partial _{k}f(x-\varepsilon
y)-\partial _{k}f(x)]dy.
\end{eqnarray*}%
Hence, by Lemma \ref{lemn1}, 
\begin{equation*}
|\partial _{kl}^{2}f^{\varepsilon }(x)|\leq C\kappa (\varepsilon ,2,\beta
)|f|_{\beta }.
\end{equation*}

To prove (\ref{maf5'}), apply (\ref{ff30}) and the interpolation theorem.
Let $\beta \in (0,1]$. Consider an operator on $C^{\beta }$ defined by $%
T^{\varepsilon }(f)=f^{\varepsilon }$. According to (\ref{ff30}), $%
T^{\varepsilon }:C^{\beta }(\mathbf{R}^{d})\rightarrow C^{k}(\mathbf{R}%
^{d}),k=1,2,$ is bounded,%
\begin{equation*}
|T^{\varepsilon }(f)|_{k}\leq C\varepsilon ^{-k+\beta }|f|_{\beta
},k=1,2,f\in C^{\beta }(\mathbf{R}^{d}).
\end{equation*}%
By Theorem 6.4.5 in \cite{BeL76}, $T^{\varepsilon }:C^{\beta }(\mathbf{R}%
^{d})\rightarrow C^{\alpha }(\mathbf{R}^{d})$ is bounded and%
\begin{equation*}
|T^{\varepsilon }(f)|_{\alpha }\leq C\varepsilon ^{(-1+\beta )(2-\alpha
)}\varepsilon ^{(-2+\beta )(\alpha -1)}|f|_{\beta }=C\varepsilon ^{-\alpha
+\beta }|f|_{\beta },f\in C^{\beta }(\mathbf{R}^{d}).
\end{equation*}%
If $\beta \in (1,\alpha ],$ $\partial ^{\alpha -1}\nabla f^{\varepsilon
}=\partial ^{\alpha -1}\left( \nabla f\right) ^{\varepsilon }$ and by (\ref%
{maf5}),%
\begin{equation*}
|\partial ^{\alpha -1}\nabla f^{\varepsilon }(x)|=|\partial ^{\alpha
-1}\left( \nabla f\right) ^{\varepsilon }(x)|\leq C\kappa (\varepsilon
,\alpha -1,\beta -1)|\nabla f|_{\beta -1},
\end{equation*}%
and (\ref{maf6}) follows.

Let $\beta \in (0,1],\alpha \in (1,2),\beta <\alpha -1.$ Then%
\begin{eqnarray*}
\nabla f^{\varepsilon } &=&\varepsilon ^{-1}\varepsilon ^{-d}\int \nabla w(%
\frac{y}{\varepsilon })f(x-y)dy \\
&=&\varepsilon ^{-1}\varepsilon ^{-d}\int \nabla w(\frac{y}{\varepsilon }%
)f(x-y)dy.
\end{eqnarray*}%
and%
\begin{equation*}
\partial ^{\alpha -1}\nabla f^{\varepsilon }=\varepsilon ^{-1-(\alpha
-1)}\int \partial ^{\alpha -1}(\nabla w)(y)f(x-\varepsilon y)dy
\end{equation*}%
and we derive as in part (i) (using Lemma \ref{rem2}) that%
\begin{equation*}
|\partial ^{\alpha -1}\nabla f^{\varepsilon }(x)|\leq C\varepsilon ^{-\alpha
+\beta }|f|_{\beta }.
\end{equation*}%
If $\beta \in (0,1],\alpha \in (1,2),\beta >\alpha -1$, then%
\begin{eqnarray*}
\partial ^{\alpha -1}\nabla f^{\varepsilon } &=&\varepsilon ^{-1}\int \nabla
w(y)(\partial ^{\alpha -1})f(x-\varepsilon y)dy \\
&=&\varepsilon ^{-1}\int \nabla w(y)[(\partial ^{\alpha -1})f(x-\varepsilon
y)-(\partial ^{\alpha -1})f(x)]dy
\end{eqnarray*}%
and%
\begin{equation*}
|\partial ^{\alpha -1}\nabla f^{\varepsilon }|\leq C\varepsilon ^{-1+\beta
-\alpha +1}|\partial ^{\alpha -1}f|_{\beta -\alpha +1}\leq C\varepsilon
^{-\alpha +\beta }|f|_{\beta }.
\end{equation*}

The statement is proved.
\end{proof}

\begin{corollary}
\label{coro3}Let $\alpha \in (0,2],\beta \leq \alpha $. Assume $\varepsilon
\in (0,1)$, $a(x)$ is bounded, and 
\begin{equation*}
\int \big(|y|^{\alpha }\wedge 1\big)\pi (dy)<\infty .
\end{equation*}%
Then there exists a constant $C$ such that for all $z,x\in \mathbf{R}^{d}$, $%
f\in C^{\beta }(\mathbf{R}^{d}),$%
\begin{equation*}
|B_{z}f^{\varepsilon }(x)|\leq C\kappa (\varepsilon ,\alpha ,\beta
)|f|_{\beta }.
\end{equation*}
\end{corollary}

\begin{proof}
If $\beta \leq \alpha <1$, by Lemmas \ref{r1} and \ref{lnew2},%
\begin{equation*}
f^{\varepsilon }(x+y)-f^{\varepsilon }(x)=\int k^{(\alpha )}(y,y^{\prime
})\partial ^{\alpha }f^{\varepsilon }(x-y^{\prime })dy^{\prime },
\end{equation*}%
and by (\ref{maf5}), 
\begin{equation*}
|f^{\varepsilon }(x+y)-f^{\varepsilon }(x)|\leq C\kappa (\varepsilon ,\alpha
,\beta )|f|_{\beta }(|y|^{\alpha }\wedge 1),x,y\in \mathbf{R}^{d}.
\end{equation*}
So,%
\begin{eqnarray*}
|f^{\varepsilon }(x+G(x)y)-f^{\varepsilon }(x)| &\leq &C\kappa (\varepsilon
,\alpha ,\beta )|f|_{\beta }(|G(x)y|^{\alpha }\wedge 1) \\
&\leq &C\kappa (\varepsilon ,\alpha ,\beta )|f|_{\beta }\big[\mathbf{1}%
_{\left\{ |y|\leq 1\right\} }|G(x)y|^{\alpha } \\
&&\quad +\mathbf{1}_{\left\{ |y|>1\right\} }(|G(x)y|^{\alpha }\wedge 1)\big].
\end{eqnarray*}

If $\beta \leq \alpha =1$, by Lemma \ref{lnew2}, (\ref{ff30}),%
\begin{eqnarray*}
|f^{\varepsilon }(x+y)-f^{\varepsilon }(x)| &\leq &C\sup_{x}[f(x)|+|\nabla
f^{\varepsilon }(x)|](|y|\wedge 1) \\
&\leq &C\kappa (\varepsilon ,1,\beta )|f|_{\beta }(|y|\wedge 1),x,y\in 
\mathbf{R}^{d}
\end{eqnarray*}%
and%
\begin{eqnarray*}
|f^{\varepsilon }(x+G(x)y)-f^{\varepsilon }(x)| &\leq &C\kappa (\varepsilon
,1,\beta )|f|_{\beta }(|G(x)y|\wedge 1) \\
&\leq &C\kappa (\varepsilon ,1,\beta )|f|_{\beta }[\mathbf{1}_{|y|\leq
1}|G(x)y|+\mathbf{1}_{|y|>1}(|G(x)y|\wedge 1)].
\end{eqnarray*}

Assume $\alpha \in (1,2],\beta \leq \alpha .$ Then for $x,y^{\prime }\in 
\mathbf{R}^{d},$%
\begin{equation}
f^{\varepsilon }(x+y^{\prime })-f^{\varepsilon }(x)-(\nabla f^{\varepsilon
}(x),y^{\prime })=\int_{0}^{1}\big(\nabla f^{\varepsilon }(x+sy^{\prime
})-\nabla f^{\varepsilon }(x),y^{\prime }\big)ds.  \label{ff34}
\end{equation}%
We will show that for all $x,y^{\prime }\in \mathbf{R}^{d}$%
\begin{eqnarray}
&&|f^{\varepsilon }(x+y^{\prime })-f^{\varepsilon }(x)-(\nabla
f^{\varepsilon }(x),y^{\prime })|  \label{ff350} \\
&\leq &C|y^{\prime }|^{\alpha }\kappa (\varepsilon ,\alpha ,\beta
)|f|_{\beta }  \notag
\end{eqnarray}%
If $\beta \in (0,\alpha ],\alpha \in (1,2),\beta \neq \alpha -1$, then we
have (\ref{ff350}) by Lemmas \ref{r1}, \ref{lnew2}$.$

If $\beta \in (0,\alpha ],\alpha \in (1,2),\beta =\alpha -1$,then for any $%
x,y^{\prime }\in \mathbf{R}^{d}$%
\begin{eqnarray*}
&&\nabla f^{\varepsilon }(x+y^{\prime })-\nabla f^{\varepsilon }(x) \\
&=&\varepsilon ^{-1}\varepsilon ^{-d}\int \nabla w(y/\varepsilon
)[f(x+y^{\prime }-y)-f(x-y)]dy
\end{eqnarray*}%
and%
\begin{equation*}
|\nabla f^{\varepsilon }(x+y^{\prime })-\nabla f^{\varepsilon }(x)|\leq
C\varepsilon ^{-1}|y^{\prime }|^{\beta }|f|_{\beta }=C\varepsilon ^{-\alpha
+\beta }|f|_{\beta }|y^{\prime }|^{\alpha -1}.
\end{equation*}%
So, (\ref{ff350}) holds in this case as well. If $\beta \leq \alpha =2$,
then by Lemma \ref{lnew2},%
\begin{eqnarray*}
|\nabla f^{\varepsilon }(x+y^{\prime })-\nabla f^{\varepsilon }(x)| &\leq
&\sup_{x}|D^{2}f^{\varepsilon }(x)|~|y^{\prime }| \\
&\leq &C\kappa (\varepsilon ,2,\beta )|f|_{\beta }|y^{\prime }|
\end{eqnarray*}%
and (\ref{ff350}) follows. Hence, for $|y|\leq 1,$ by (\ref{ff350}),%
\begin{equation*}
|f^{\varepsilon }(x+G(x)y)-f^{\varepsilon }(x)-(\nabla f^{\varepsilon
}(x),G(x)y)|\leq C\kappa (\varepsilon ,\alpha ,\beta )|G(x)y|^{\alpha
}|f|_{\beta }.
\end{equation*}%
Also, for $|y|>1,$%
\begin{equation*}
|f^{\varepsilon }(x+G(x)y)-f^{\varepsilon }(x)|\leq 2|f|_{\beta }.
\end{equation*}

Therefore, the statement follows by the assumptions and Lemma \ref{lnew2}.
\end{proof}


\subsection{Proof of Lemma \protect\ref{lem:expect}}

If $\beta \leq \alpha \,,$ define $f^{\varepsilon }$ by (\ref{maf7}) for $%
\varepsilon \in (0,1)$ and apply It\^{o}'s formula (see Remark \ref{lrenew1}%
): for $s\in \lbrack 0,T]$, 
\begin{equation*}
\mathbf{E}[f^{\varepsilon }(Y_{s})-f^{\varepsilon }(Y_{\tau _{i_{s}}})|%
\mathcal{F}_{\tau _{i_{s}}}]=\mathbf{E}\big[\int_{\tau _{i_{s}}}^{s}\big(%
A_{Y_{\tau _{i_{s}}}}f^{\varepsilon }(Y_{r})+B_{Y_{\tau
_{i_{s}}}}f^{\varepsilon }(Y_{r})\big)dr\big|\mathcal{F}_{\tau _{i_{s}}}\big]%
.
\end{equation*}%
Hence, by Lemma \ref{lnew2} and Corollary \ref{coro3}, for $\varepsilon \in
(0,1)$,%
\begin{eqnarray*}
|\mathbf{E}[f(Y_{s})-f(Y_{\tau _{i_{s}}})|\mathcal{F}_{\tau _{i_{s}}}]|
&\leq &|\mathbf{E}[(f-f^{\varepsilon })(Y_{s})-(f-f^{\varepsilon })(Y_{\tau
_{i_{s}}})|\mathcal{F}_{\tau _{i_{s}}}]| \\
&&+|\mathbf{E}[f^{\varepsilon }(Y_{s})-f^{\varepsilon }(Y_{\tau _{i_{s}}})|%
\mathcal{F}_{\tau _{i_{s}}}]| \\
&\leq &CF(\varepsilon ,\delta )|f|_{\beta },
\end{eqnarray*}%
with a constant $C$ independent of $\varepsilon ,f$ and 
\begin{equation*}
F(\varepsilon ,\delta )=\left\{ 
\begin{array}{cc}
(\varepsilon ^{2}+\delta )(1-\ln \varepsilon ) & \text{if }\alpha =\beta =2,
\\ 
\varepsilon ^{\beta }+\delta \kappa (\varepsilon ,\alpha ,\beta ) & \text{%
otherwise.}%
\end{array}%
\right.
\end{equation*}%
Minimizing $F(\varepsilon ,\delta )$ in $\varepsilon \in (0,1)$, we obtain%
\begin{equation*}
|\mathbf{E}[f(Y_{s})-f(Y_{\tau _{i_{s}}})|\mathcal{F}_{\tau _{i_{s}}}]|\leq
Cr(\delta ,\alpha ,\beta )|f|_{\beta }.
\end{equation*}

If $\beta >\alpha ,$ apply It\^{o}'s formula directly (see Remark \ref%
{lrenew1}):%
\begin{equation*}
\mathbf{E}[f(Y_{s})-f(Y_{\tau _{i_{s}}})|\mathcal{F}_{\tau _{i_{s}}}]=%
\mathbf{E}\big[\int_{\tau _{i_{s}}}^{s}\big(A_{Y_{\tau _{i_{s}}}}^{(\alpha
)}f(Y_{r})+B_{Y_{\tau _{i_{s}}}}^{(\alpha )}f(Y_{r})\big)dr \big|\mathcal{F}%
_{\tau _{i_{s}}}\big].
\end{equation*}%
Hence, by Lemmas \ref{prop2} and \ref{lnew2}, 
\begin{equation*}
|\mathbf{E}[f(Y_{s})-f(Y_{\tau _{i_{s}}})|\mathcal{F}_{\tau _{i_{s}}}]|\leq
C\delta |f|_{\beta }.
\end{equation*}%
The statement of Lemma \ref{lem:expect} follows.


\subsection{Proof of Theorem \protect\ref{thm:main}}

Let $v\in C^{\alpha +\beta }(H)$ be the unique solution to (\ref{maf8}) (see
Corollary \ref{lcornew1}). By It\^{o}'s formula (see Remark \ref{lrenew1})
and (\ref{maf8}),%
\begin{eqnarray*}
\mathbf{E}[v(0,X_{0})] &=&\mathbf{E}[v(T,X_{T})] - \mathbf{E}\big[%
\int_{0}^{T}\big(\partial
_{t}v(s,X_{s})+A_{X_{s}}v(s,X_{s})+B_{X_{s}}v(s,X_{s})\big)ds\big] \\
&=&\mathbf{E}\big[g(X_{T})-\int_{0}^{T}f(X_{s})ds\big]
\end{eqnarray*}%
and 
\begin{equation}
\mathbf{E}[v(0,X_{0})] = \mathbf{E}[v(0,Y_{0})].  \label{eqn:expect_terminal}
\end{equation}

By Proposition \ref{b1}, Corollary~\ref{lcornew1}, Remark \ref{rlast}, and
Lemma \ref{prop2}, 
\begin{eqnarray}
|A_{z}v(s,\cdot )|_{\beta }+|B_{z}v(s,\cdot )|_{\beta } &\leq &C|v|_{\alpha
+\beta }\leq C|g|_{\alpha +\beta },  \label{maf9} \\
|\partial _{t}v(s,\cdot )|_{\beta } &\leq &C|g|_{\alpha +\beta },s\in
\lbrack 0,T].  \notag
\end{eqnarray}

Then, by It\^{o}'s formula (Remark \ref{lrenew1}) and Corollary \ref%
{lcornew1}, with (\ref{eqn:expect_terminal}) and (\ref{maf9}), it follows
that 
\begin{eqnarray*}
&&\mathbf{E}[g(Y_{T})] - \mathbf{E}[g(X_{T})] - \mathbf{E}\big[%
\int_{0}^{T}f(Y_{\tau _{i_{s}}})ds\big] + \mathbf{E}\big[%
\int_{0}^{T}f(X_{s})ds\big] \\
&=&\mathbf{E}[v(T,Y_{T})] - \mathbf{E}[v(0,Y_{0})] - \mathbf{E}\big[%
\int_{0}^{T}f(Y_{_{\tau _{i_{s}}}})ds\big] + \mathbf{E}\big[%
\int_{0}^{T}f(X_{s})ds\big] \\
&=&\mathbf{E}\Big[\int_{0}^{T}\Big\{\big[\partial _{t}v(s,Y_{s})-\partial
_{t}v(s,Y_{\tau _{i_{s}}})\big] \\
&&+\big[A_{Y_{\tau _{i_{s}}}}v(s,Y_{s})-A_{Y_{\tau _{i_{s}}}}v(s,Y_{\tau
_{i_{s}}})\big] \\
&&+\big[B_{Y_{\tau _{i_{s}}}}v(s,Y_{s})-B_{Y_{\tau _{i_{s}}}}v(s,Y_{\tau
_{i_{s}}})\big]\Big\}ds\Big].
\end{eqnarray*}

Hence, by (\ref{maf9}) and Lemma~\ref{lem:expect}, there exists a constant $%
C $ independent of $g$ such that 
\begin{equation*}
|\mathbf{E}[g(Y_{T})]-\mathbf{E}[g(X_{T})]|\leq Cr(\delta ,\alpha ,\beta
)|g|_{\alpha +\beta }.
\end{equation*}%
The statement of Theorem \ref{thm:main} follows.

\section{Conclusion}

The paper studies weak Euler approximation of SDEs driven by L\'{e}vy
processes. The dependence of the rate of convergence on the regularity of
coefficients and driving processes is investigated under assumption of $%
\beta $-H\"{o}lder continuity of the coefficients. It is assumed that the
main term of the SDE is driven by a spherically-symmetric $\alpha $-stable
process and the tail of the L\'{e}vy measure of the lower order term has a $%
\mu $-order finite moment ($\mu \in (0,3)).$ The resulting rate depends on $%
\beta ,\alpha $ and $\mu $. In order to estimate the rate of convergence,
the existence of a unique solution to the corresponding backward Kolmogorov
equation in H\"{o}lder space is first proved. The assumptions on the
regularity of coefficients and test functions are different than those in
the existing literature.

One possible improvement could be to consider the asymptotics of the tails
at infinity instead of the tail moment $\mu $. Besides this, the stochastic
differential equations considered so far are associated with nondegenerate L%
\'{e}vy operators. A further step could be to study the case with degenerate
operators. That is, consider equation (\ref{one}) without assuming $\det
b\neq 0$. For example, let $\alpha \in \lbrack 1,2]$ and $\beta \in (\alpha
,2\alpha ].$ Assume the coefficients are in $C^{\beta }$ and 
\begin{equation*}
\int_{|y|\leq 1}|y|^{\alpha }d\pi +\int_{|y|>1}|y|^{2\alpha }d\pi <\infty .
\end{equation*}%
In this case, a plausible convergence rate is $r(\delta ,\alpha ,\beta
)=\delta ^{\frac{\beta }{\alpha }-1}$. With $\det b=0$ being allowed, a
higher regularity of coefficients and lighter tails of $\pi $ would be
required.


\end{document}